\documentclass[letterpaper, 10pt, journal, twocolumn, final]{IEEEtran}

\usepackage{cite}
\usepackage{algorithm}
\usepackage{algpseudocode}
\usepackage{algorithmicx}
\usepackage{graphicx}
\usepackage{subfig}
\usepackage{textcomp}

\usepackage{amssymb,amsmath,amsfonts,cases,mathtools}
\usepackage{microtype}
\usepackage{url}
\usepackage[hidelinks]{hyperref}
\usepackage{xifthen}
\usepackage{xpatch}
\usepackage{amsthm}
\usepackage[dvipsnames]{xcolor}
\usepackage{tikz}
\usepackage{lipsum,booktabs}
\usepackage{pifont}
\usepackage{xcolor}
\usepackage{soul}
\usepackage{balance}

\newcommand{\R}{\mathbb{R}}

\newcommand{\N}{\mathbb{N}}

\newcommand{\T}{^\top}
\newcommand{\inv}{^{-1}}

\newcommand{\subj}{\tn{subj.~to}}

\DeclareMathOperator{\col}{col}

\DeclareMathOperator{\Tr}{Tr}

\newtheorem{theorem}{Theorem}

\newtheorem{lemma}{Lemma}
\newtheorem{assumption}{Assumption}
\newtheorem{remark}{Remark}

\usepackage[dvipsnames]{xcolor}
\usepackage{tikz,pgfplots}
\pgfplotsset{compat=newest}

\newcommand\oprocendsymbol{\hbox{$\blacksquare$}}
\newcommand\oprocend{\relax\ifmmode\else\unskip\hfill\fi\oprocendsymbol}

\newcommand{\Tsim}{T}

\newcommand{\tn}[1]{{\textnormal{#1}}}

\newcommand{\cB}{\mathcal B}
\newcommand{\cC}{\mathcal C}

\newcommand{\cK}{\mathcal K}

\newcommand{\cQ}{\mathcal Q}

\newcommand{\cS}{\mathcal S}

\newcommand{\cV}{\mathcal V}

\newcommand{\cX}{\mathcal X}

\newcommand{\until}[1]{\{1,\ldots,#1\}} 

\graphicspath{{figs/}}

\newcommand{\norm}[1]{\left \|#1 \right \|}

\newcommand{\rz}{\mathrm{z}}

\newcommand{\uu}{u}
\newcommand{\z}{\rz}

\newcommand{\dd}{\mathrm{D}}

\newcommand{\ud}{^}%
\newcommand{\du}{_}
\newcommand{\iter}{k}

\newcommand{\initer}{\tau}
\newcommand{\iterp}{{\iter+1}}
\newcommand{\timeid}{t}
\newcommand{\timeidp}{{\timeid + 1}}

\newcommand{\utime}{\uu\du{\timeid}}

\newcommand{\av}{_{\textsc{av}}}

\newcommand{\x}{x}

\newcommand{\xtime}{\x\du\timeid}
\newcommand{\xtimep}{\x\du{\timeidp}}

\newcommand{\dr}{\dd}
\newcommand{\drt}{\dr\ud\iter}

\newcommand{\g}{K} %
\newcommand{\ga}{\g\av} 

\newcommand{\gt}{\g\ud{\iter}}
\newcommand{\gtp}{\g\ud{\iterp}}
\newcommand{\gat}{\ga\ud\iter}
\newcommand{\gatp}{\ga\ud{\iterp}}

\newcommand{\gstar}{\g^\star}

\newcommand{\G}{G}
\newcommand{\J}{J}

\newcommand{\Ja}{J\av^\delta}

\newcommand{\lipp}{\beta}

\newcommand{\levset}{\Omega}
\newcommand{\tlevset}{\tilde{\levset}}
\newcommand{\levsetg}{\levset^{\g}}

\newcommand{\dimx}{{n}}
\newcommand{\dimu}{{m}}
\newcommand{\bE}{\mathbb{E}}

\newcommand{\gd}{\mu}
\newcommand{\gl}{\psi}

\newcommand{\mav}{\lambda}

\newcommand{\ssz}{\gamma}

\newcommand{\f}{\mathrm{z}} %
\newcommand{\ft}{\f\ud\iter}
\newcommand{\ftp}{\f\ud\iterp}

\newcommand{\fa}{\f\av} 

\newcommand{\fat}{\fa\ud\iter} 
\newcommand{\fatp}{\fa\ud\iterp}

\newcommand{\tfa}{\tilde{\f}\av} 

\newcommand{\tfat}{\tfa\ud\iter}
\newcommand{\tfatp}{\tfa\ud\iterp}

\newcommand{\pert}{g^\delta}

\newcommand{\Jt}{\J_{\Tsim}}

\newcommand{\acc}{\alpha}

\newcommand{\app}{p^{\delta,\Tsim}}

\newcommand{\lippa}{\lipp_{p}}

\newcommand{\lippl}{\lipp_{\err}}

\newcommand{\tV}{\tilde{V}}

\newcommand{\lippp}{\lipp_{\J\av}}

\newcommand{\Tper}{\mathbf{\iter}_{\text{prd}}}
\newcommand{\Tpper}{\mathbf{\iter}_{p,\text{prd}}}
\newcommand{\Tpqer}{\mathbf{\iter}_{q,\text{prd}}}
\newcommand{\Trper}{\mathbf{\iter}_{r,\text{prd}}}
\newcommand{\Toper}{\mathbf{\iter}_{11,\text{prd}}}
\newcommand{\Tnmper}{\mathbf{\iter}_{\dimu\dimx,\text{prd}}}
\newcommand{\Tijper}{\mathbf{\iter}_{ij,\text{prd}}}

\newcommand{\ts}{\Delta t}

\newcommand{\Ac}{A_{\text{cont}}}
\newcommand{\Bc}{B_{\text{cont}}}

\newcommand{\sszcode}{10^{-7}} 
\newcommand{\deltacode}{10^{-2}}
\newcommand{\Tpercode}{19} 
\newcommand{\Tsimcode}{20}
\newcommand{\tscode}{10^{-2}}
\newcommand{\Qcode}{2}

\newcommand{\itercodea}{1} 
\newcommand{\itercodeb}{100} 
\newcommand{\itercodec}{1000} 
\newcommand{\itercoded}{2\cdot10^{5}}

\newcommand{\rav}{\rho\av}
\newcommand{\err}{e}

\newcommand{\vf}{\cV}
\newcommand{\vft}{\cV\du{\timeid}}
\newcommand{\vftp}{\cV\du{\timeidp}}

\newcommand{\vfT}{\vf\du{\Tsim}}

\def\algoext/{EXtremum-seeking Policy iteration LQR}

\def\algo/{EXP-LQR}

\newcommand{\state}[1]{\chi\ud{#1}}
\newcommand{\stateinit}{\chi_0}

\newcommand{\V}{V_{\mav}}

\usepackage{soul}

\newcommand{\mathsout}[1]%
{\bgroup\mathchoice
  {\sbox0{$\displaystyle{#1}$}%
    \usebox0\hspace{-\wd0}%
    \rule[0.5\ht0-0.5\dp0-.5pt]{\wd0}{1pt}}%
  {\sbox0{$\textstyle{#1}$}%
    \usebox0\hspace{-\wd0}%
    \rule[0.5\ht0-0.5\dp0-.5pt]{\wd0}{1pt}}%
  {\sbox0{$\scriptstyle{#1}$}%
    \usebox0\hspace{-\wd0}%
    \rule[0.5\ht0-0.5\dp0-.5pt]{\wd0}{1pt}}%
  {\sbox0{$\scriptscriptstyle{#1}$}%
    \usebox0\hspace{-\wd0}%
    \rule[0.5\ht0-0.5\dp0-.5pt]{\wd0}{1pt}}%
\egroup}

\usepackage[textwidth=1.4cm,textsize=scriptsize]{todonotes}
\title{Data-Driven LQR with Finite-Time Experiments\\
via Extremum-Seeking Policy Iteration
}

\author{Guido Carnevale,
Nicola Mimmo,
Giuseppe Notarstefano 
\thanks{Work supported in part by Fondi PNRR - Bando PE - Progetto
PE11 - 3A-ITALY, “Made in Italy Circolare e Sostenibile” - Codice
PE0000004, CUP: J33C22002950001 and by MOST - Sustainable Mobility National Research Center and received funding from the European Union Next-GenerationEU (PIANO NAZIONALE DI RIPRESA E RESILIENZA (PNRR) - MISSIONE 4 COMPONENTE 2, INVESTIMENTO 1.4 - D.D. 1033 17/06/2022, CN00000023). 
The authors are with the Department of Electrical,  Electronic and Information Engineering,  Alma Mater Studiorum - Universita` di Bologna,  Bologna, Italy, e-mail: 
\{guido.carnevale, nicola.mimmo2, giuseppe.notarstefano\}@unibo.it.%
}
}

\begin{document}
	\maketitle

\begin{abstract}
	In this paper, we address Linear Quadratic Regulator (LQR) problems through a novel iterative algorithm named \algoext/ (\algo/). 
	The peculiarity of \algo/ is that it only needs access to a truncated approximation of the infinite-horizon cost associated to a given policy. 
	Hence, \algo/ does not need the direct knowledge of neither the system and cost matrices. %
	In particular, at each iteration, \algo/ refines the maintained policy using a truncated LQR cost retrieved by performing finite-time virtual or real experiments in which a perturbed version of the current policy is employed.
	Such a perturbation is done according to an extremum-seeking mechanism and makes the overall algorithm a time-varying nonlinear system.
	By using a Lyapunov-based approach exploiting averaging theory, we show that \algo/ exponentially converges to an arbitrarily small neighborhood of the optimal gain matrix.
	We corroborate the theoretical results with numerical simulations involving the control of an induction motor.
\end{abstract}

\section{Introduction}

Data-driven strategies for optimal control have become an increasingly prominent trend in recent years, see, e.g., the survey~\cite{hou2013model}.
The distinctive feature of these methods stands in refining the control policy by gathering data rather than using a priori knowledge of the system. 
A key distinction in this field is between off-policy methods, where the tentative policy is not concurrently applied to the system, and on-policy methods, where the policy is implemented.

A branch of off-policy methodologies originated by the so-called Kleinman algorithm~\cite{kleinman1968on}, see, e.g., the related works~\cite{qin2014online,modares2016optimal,pang2018data,krauth2019finite,pang2021robust,lopez2023efficient}.  
We can further classify off-policy methods by distinguishing between indirect approaches~\cite{dean2019safely,mania2019certainty,ferizbegovic2019learning}, which incorporate an initial identification step before the policy formulation, and direct approaches, where data is directly applied during the policy design~\cite{de2019formulas,van2020data, rotulo2020data}. 
Direct methods have been also extended to deal with unknown linear systems with switching time-varying dynamics~\cite{rotulo2022online}, noisy data~\cite{van2020noisy,de2021low,dorfler2023certainty}, and robustness issues~\cite{berberich2020robust}.
The works~\cite{formentin2018core,iannelli2020structured,dorfler2022bridging} try to bridge the gap between indirect and direct paradigms.
Policy-gradient methods are another widely-used class of strategies, whose distinctive feature consists of optimizing the control policies through gradient-based updates, see the works~\cite{fazel2018global,bu2019lqr,zhang2020policy,zhang2020global,zhang2023revisiting,hu2023toward}.
As for the on-policy approaches, we mention the works~\cite{vrabie2009adaptive,jiang2012computational,possieri2022value}.
Recently, on-policy methods using adaptive control tools have been provided in~\cite{borghesi2023on,borghesi2024mr}.
While, in~\cite{sforni2023on,sforni2024stability,song2024role}, on-policy strategies are obtained including learning mechanisms based on the recursive least squares mechanism.
As we will detail later, our approach is based on the so-called extremum-seeking mechanism, see, e.g., the recent survey~\cite{scheinker2024100} and the works~\cite{wittenmark1995adaptive, teel2001solving, ariyur2003real, krstic2000stability, tan2006non}.
In the context of linear optimal control, extremum-seeking has been already used in~\cite{frihauf2013finite}, where, however, it is employed with the goal of finding a sequence of open-loop control steps minimizing a finite-time horizon problem.

\begin{figure}[htpb]
	\centering
	\includegraphics{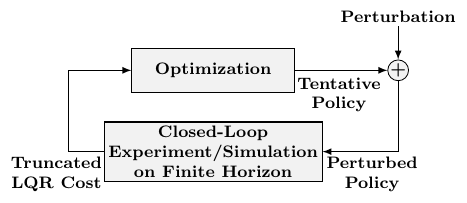}
	\caption{Schematic representation of the proposed strategy.}
	\label{fig:sketch}
\end{figure}

The main contribution of this paper is the development of \algoext/ (\algo/), namely, a novel data-driven strategy for solving LQR problems. %
Our approach does not need direct knowledge of system and cost matrices.
More specifically, our method only needs a finite-time truncated version of the infinite-horizon cost (obtained, e.g., by running the real system or a simulator) computed by using a suitably perturbed version of the current policy maintained by the algorithm, see the schematic representation provided in Fig.~\ref{fig:sketch}.
Using this information, \algo/ iteratively improves the policy taking on an extremum-seeking mechanism and a suitable reformulation of the LQR problem.
Building on the extremum-seeking framework, our mechanism employs deterministic dither matrices to perturb the current gain and leverages a low-pass filter to enhance the gradient approximation.
The design of the dither matrix, along with the incorporation of the low-pass filter and an ad-hoc policy evaluation procedure, represents a key novelty and advantage of our method compared to existing schemes in the literature, which are mostly based on random perturbations and do not exploit any filtering.
We interpret the overall algorithm as a nonlinear time-varying system, which we then analyze by using system-theoretic tools based on the so-called averaging approach (see, for example,~\cite[Ch.~10]{khalil2002nonlinear} and~\cite{sanders2007averaging} for the continuous-time case or~\cite{bai1988averaging} for the discrete-time one).
Indeed, as customary in the context of averaging theory, we focus on the so-called \emph{averaged system} associated to the algorithm.
In particular, the averaged system reads as a policy gradient method perturbed by errors arising from the use of the truncated cost instead of the infinite-horizon one, as well as from the derivative-free gradient approximation.
More in detail, we employ a Lyapunov-based approach to ensure that the averaged system trajectories exponentially converge to an arbitrarily small neighborhood of the optimal gain matrix.
Then, we use this preparatory result to achieve the same property on the trajectories of the original time-varying algorithm.
This final step is supported by Theorem~\ref{th:averaging}, introduced in Section~\ref{sec:averaging_theory}, which presents averaging-related stability results for generic discrete-time systems. 
To the best of the authors' knowledge, Theorem~\ref{th:averaging} also represents a per se contribution of this work. 
A conference version of this paper appeared in~\cite{carnevale2024extremum}. 
However, in that preliminary version, the algorithm relies on oracles providing the exact infinite-horizon cost associated to the tentative gain, making it impractical for real-world scenarios where only finite-time virtual or real experiments are feasible.
Moreover, certain proofs were omitted.
Finally, this work includes a concrete application example involving the control of an induction motor.

We organize the paper as follows.
In Section~\ref{sec:averaging_theory}, we introduce some preliminaries about averaging theory for discrete-time systems.
In Section~\ref{sec:problem_setup}, we describe the problem setup considered in the paper.
In Section~\ref{sec:algo}, we provide the description of \algo/ and state its theoretical features.
Finally, in Section~\ref{sec:analysis}, we numerically test the effectiveness of \algo/.

\paragraph*{Notation}
A square matrix $M \in \R^{n \times n}$ is Schur if all its eigenvalues lie in the open unit disk.
The identity matrix in $\R^{n \times n}$ is $I_n$. 
The vector of zeros of dimension $n$ is denoted as $0_n$.
The vertical concatenation of vectors $v_1, \dots, v_N$ is $\col(v_1, \dots, v_N)$.
Given $r > 0$ and $x \in \R^n$, we use $\cB_r(x)$ to denote the closed ball of radius $r > 0$ centered in $x$, namely $\cB_r(x) := \{y \in \R^n \mid \norm{y - x} \leq r\}$.
Given $A \in \R^{n \times n}$, $\Tr(A)$ denotes its trace.
$\R_+$ denotes the positive orthant in $\R$.

\section{Preliminaries: Averaging Theory for Discrete-Time Systems}
\label{sec:averaging_theory}

In this preliminary part, we provide a generic stability result for
discrete-time systems in the context of averaging theory (see, e.g.,~\cite{khalil2002nonlinear,sanders2007averaging,bai1988averaging}).
Although we will use it as an instrumental step for proving the main
result of the paper, we remark that it represents a contribution  per se.

Let us consider the time-varying discrete-time system 
\begin{equation}
    \label{eq:plant}
    \state\iterp = \state\iter +\ssz f(\state\iter ,\iter) \qquad \state0 = \stateinit,
\end{equation}
where $\state\iter \in \R^n$ denotes the state, $f: \R^{n} \times \N \to \R^n$ describes its dynamics, and $\ssz > 0$ is a tunable parameter.  
Let us enforce the following assumptions. 
\begin{assumption}
\label{hyp:fa}
    There exist $\Tper  \in \N$ and  $f\av\,:\,\R^n\to\R^n$ such that
\begin{equation}
\label{eq:fa}
    f\av (\state{}) = \dfrac{1}{\Tper } \sum_{\initer = \iter+1}^{\iter+\Tper} f(\state{},\initer),
\end{equation}
for all $\state{} \in \R^n$ and $\iter\in\N$.\oprocend
\end{assumption}
Assumption~\ref{hyp:fa} allows for properly writing a well-posed averaged system associated to system~\eqref{eq:plant}.
Roughly, Assumption~\ref{hyp:fa} says that $f(\chi,\cdot)$ is periodic and $\Tper \in \N$ represents its period.
The next assumption guarantees some regularity conditions on the functions $f$ and $f\av$ and their derivatives.
\begin{assumption}\label{hyp:bounds}
  There exists a set $\cX \subseteq \R^n$ such that the restrictions of $f(\cdot,k)$, $f\av(\cdot)$, $\partial f(\chi,k)/\partial \chi$, and $\partial f\av(\chi) / \partial \chi$ to $\mathcal{X}$ are continuous for all $k \in \N$.
  \oprocend
\end{assumption}
The next assumption characterizes the convergence properties of the so-called \emph{averaged system} associated to~\eqref{eq:plant}, i.e., the auxiliary time-invariant dynamics of $\state\iter\av \in \R^{n}$ described by 
\begin{equation}\label{eq:average}
	\state\iterp\av = \state\iter\av +\ssz f\av(\state\iter\av) \qquad \state0\av = \stateinit.
\end{equation}
To this end, we first introduce a continuously differentiable function $V: \cX \to \R_+$ and, given any $c > 0$, its level set defined as %
$\Omega_c := \{x \in \R^n \mid V(x) \leq c\}$.
\begin{assumption}
\label{hyp:StabilityAv}
	For all $c_0 > 0$ and all $\rav \in (0, c_0)$, there exist $\bar{\gamma}_1 > 0$ and $a \in (0,1)$ such that, for all $\chi_0 \in \Omega_{c_0}$ and $\gamma \in (0,\bar{\gamma_1})$, it holds
	\begin{equation}\label{eq:StabilityAv}
		V(\state{\iter}\av) \leq (1-\ssz a)^\iter V(\stateinit) + \rav,
	\end{equation}    
	along the trajectories of~\eqref{eq:average} for all $\iter \in \N$.
	\oprocend
\end{assumption}
We are ready to state the following result about the original system~\eqref{eq:plant}.
\begin{theorem}\label{th:averaging}
	Consider system~\eqref{eq:plant} and let Assumptions~\ref{hyp:fa}-\ref{hyp:StabilityAv} hold. 
	Then, for all $c_0, c_1 > 0$ such that $\Omega_{c_0} \subset \Omega_{c_1} \subseteq \mathcal{X}$ and $\rho \in (0, c_0)$, there exists $\bar{\ssz} > 0$ such that, for all $\chi_0 \in \Omega_{c_0}$ and $\ssz \in (0,\bar{\ssz})$, it holds
	\begin{subequations}
		\begin{align}
			V(\state\iter)   &\leq c_1,\label{eq:bounded_trajectory}
			\\
			V(\state{\iter}) &\le (1-\ssz a)^\iter V(\stateinit) + \rho,\label{eq:final_radius}
		\end{align}
	\end{subequations}
	for all $\iter\in\N$.
	\oprocend
\end{theorem}
The proof of Theorem~\ref{th:averaging} is provided in Appendix~\ref{sec:proof_averaging_generic}.
Essentially, Theorem~\ref{th:averaging} ensures that, with sufficiently small values of the parameter $\ssz$, the properties of the averaged system~\eqref{eq:average} enforced by Assumption~\ref{hyp:StabilityAv} can be ``transferred'' to the original time-varying system~\eqref{eq:plant}.

\section{Problem Setup}
\label{sec:problem_setup}

This section states the problem setup that we aim to address and recalls a model-based iterative approach to solve it.
	
\subsection{Data-Driven LQR Problem Setup}
\label{sec:data_driven_LQR}

In this paper, we focus on LQR problems in the form
\begin{subequations}	
	\label{eq:problem}
	\begin{align}
		\min_{\substack{\x_1, \x_2,\ldots,\\ \uu_0,\uu_1,\ldots}} \: \: & \: \bE\bigg[\frac{1}{2}\sum_{\timeid=0}^{\infty} \Big( {\xtime}\T Q \xtime + {\utime}\T R \utime  \Big)\bigg]
		\\
		\subj \: & \: \xtimep = A \xtime + B \utime, \qquad \qquad \x\du{0} \sim \cX\du{0},
		\label{eq:dyn}
	\end{align}
\end{subequations}	
where  $\xtime \in \R^\dimx$ and $\utime \in \R^\dimu$ denote, respectively, the state and the input of the system at time $\timeid \in \N$, $A \in \R^{\dimx \times \dimx}$ and $B \in \R^{\dimx \times \dimu}$ represent the state and the input matrices, while $Q \in \R^{\dimx \times \dimx}$ and $R \in \R^{\dimu \times \dimu}$ are the cost matrices.
As for the initial condition $\x\du{0} \in \R^\dimx$, we assume that it is drawn from the uniform probability distribution $\cX\du{0}$ over the unitary-sphere. 
The operator $\bE[\cdot]$ denotes the expected value with respect to
$\cX\du{0}$. %
We require the following properties on the pairs $(A,B)$ and $(Q,R)$.
\begin{assumption}[System and Cost Matrices Properties]
	\label{ass:A_B_Q_R}
	The pair $(A,B)$ is controllable, while the cost matrices $Q$
        and $R$ are both symmetric and positive definite, i.e., $Q =
        Q\T \succ 0$ and $R = R\T \succ 0$.\oprocend 
\end{assumption}
Under the properties enforced by Assumption~\ref{ass:A_B_Q_R}, when $(A,B)$ and $(Q,R)$ are known, the optimal solution to problem~\eqref{eq:problem} is ruled by a linear time-invariant policy $\utime = \gstar \xtime$ with $\gstar \in \R^{\dimu \times \dimx}$ given by
\begin{align*}
	\g^\star = - (R + B\T P^\star B)\inv B\T P^\star A,
\end{align*}
where the matrix $P^\star \in \R^{\dimx \times \dimx}$ solves the so-called Discrete-time Algebraic Riccati Equation associated to problem~\eqref{eq:problem}, see~\cite{anderson2007optimal}.
However, as formalized in the next assumption, in this paper the knowledge of the pairs $(A,B)$ and $(Q,R)$ is not available and, therefore, $\g^\star$ cannot be computed.
\begin{assumption}[Unknown System and Cost Matrices]
	\label{ass:unknown}
	The pairs $(A,B)$ and $(Q,R)$ are unknown.
	\oprocend
\end{assumption}
Accordingly, we are interested in devising a data-driven strategy to
iteratively address problem~\eqref{eq:problem}.

\subsection{Model-based Gradient Method for LQR}

Next, we recall a model-based gradient method to address problem~\eqref{eq:problem} in an iterative fashion.
Let $\cK \subset \R^{m \times n}$ be the set of stabilizing gains, namely 
\begin{align*}
	\cK := \{\g \in \R^{m \times n} \mid A + B\g \text{ is Schur}\}.
\end{align*}
As shown in, e.g.,~\cite{bu2019lqr}, by considering the state-feedback control $\utime = \g\xtime$ with $\g \in \cK$, it is possible to recast problem~\eqref{eq:problem} as the unconstrained program
\begin{align}\label{eq:LQR}
	\min_{\g \in \cK}%
	\: & \: \J(\g),
\end{align} 
where the cost function $\J: \cK  \to \R$ is given by
\begin{align}
	& \J(\g) \! := \! 
	\frac{1}{2\dimx}\!
	\Tr\!\bigg(\!
	\sum_{\timeid=0}^\infty
	(A+B\g)^{\timeid,\top}
	(Q \! + \! \g\T R \g)
	(A+B\g)^{\timeid}\! \bigg)\!.\notag %
\end{align}
It is worth noting that since $\x\du{0} \sim \cX\du{0}$ (see problem~\eqref{eq:problem}) and $\cX\du 0$ is a uniform distribution over the unitary-radius sphere,
then the set of stabilizing gains $\cK$ coincides with the domain of the cost function $\J$~\cite{bu2019lqr}.
Moreover, being the set $\cK$ open~\cite[Lemma~IV.3]{bu2020topological} and connected~\cite[Lemma~IV.6]{bu2020topological}, one could use the gradient descent method to iteratively solve problem~\eqref{eq:LQR} (see, e.g.,~\cite{bu2019lqr}).
Namely, at each iteration $\iter \in \N$, an estimate $\gt \in \R^{\dimx \times \dimu}$ of the optimal gain $\gstar$ could be maintained and iteratively updated according to
\begin{align}
	\gtp = \gt - \ssz \G(\gt),
	\label{eq:desired_update}
\end{align}
where $\ssz > 0$ is the step size parameter, while, when $\R^{m \times n}$ is equipped with the Frobenius inner product, $\G: \R^{m \times n} \to \R^{m \times n}$ is the gradient of the cost function $J$ with respect to $\g$ evaluated at $\gt$.
In particular, given $\g \in \cK$, the gradient $\G(\g)$ reads as
\begin{align*}
	&\G(\g) = \left( R \g +B\T P(A+B\g) \right) W_c,
\end{align*}
where the matrices $W_c\in \R^{n \times n}$ and $P \in \R^{n \times n}$ are solutions to the equations
\begin{align*}
		(A+B\g) W_c (A + B\g)\T - W_c&=-I_n
		\\
		(A+B\g)\T P (A+B\g) - P &= -(Q+{\g}^\top R \g).
\end{align*}
Hence, in our setup, it is not possible to compute $\G(\gt)$ and implement~\eqref{eq:desired_update} because its computation would require the knowledge of the pairs $(A,B)$ and $(Q,R)$ that are both not available (cf. Assumption~\ref{ass:unknown}).
However, for a given gain $\g$ (e.g., the current estimate about the optimal gain $\gstar$), we assume the presence of an oracle providing the finite-horizon cost
\begin{align*}
	\Jt(\g) \! := \! 
\frac{1}{2\dimx}\!\Tr\!\left(\sum_{\timeid=0}^{\Tsim-1}
(A \! + \! B\g)^{\timeid,\top}
(Q \! + \! \g\T R \g)
(A \! + \! B\g)^{\timeid}\!\right)\!\!,
\end{align*}
where the number of samples $\Tsim \in \N$ represents an algorithm parameter that will be designed later.
Differently from the entire cost $\J(\g)$ whose exact computation would require virtual or real experiments over infinite-time horizons, we remark that $\Jt(\g)$ may be retrieved with finite-time virtual or real experiments using the control law $\utime = \g\xtime$.
Since the initial condition $\x\du0$ is drawn from the uniform distribution $\cX\du0$ over the unitary sphere, $\Jt(\g)$ can be \emph{exactly} computed as the mean truncated cost achieved over $\dimx$ experiments, namely
\begin{align}
	\Jt(\g) = \frac{1}{\dimx}\sum_{i=1}^{\dimx} \vfT(\g,e_i),\label{eq:exact_cost}
\end{align}
where 
$\vfT(\g,e_i)$ is the truncated LQR cost obtained by running system~\eqref{eq:dyn} with control gain $\g$ and initial condition $e_i$, i.e., the $i$-th canonical basis vector, namely
\begin{subequations}\label{eq:cost_routine}
	\begin{align}
		\vftp(\g,e_i) &= \vft(\g,e_i) + \frac{1}{2}\xtime\T\left(Q + \g\T R\g\right)\xtime
		\\
		\xtimep &= (A + B\g)\xtime,
	\end{align}
\end{subequations}
for all $\timeid \in \{0,\dots,\Tsim-1\}$ with $\vf\du0(\g,e_i) = 0$ and $\x\du0 = e_i$.
\begin{remark}\label{rem:Q_R}
	Notice that the availability of the experimental truncated costs $\vfT(\g,e\du{1}), \dots, \vfT(\g,e\du{\dimx})$ does not necessarily imply the knowledge of $(Q, R)$.
	For instance, the experimental costs $\vfT(\g,e_i)$ in~\eqref{eq:cost_routine} may be obtained from sensor data or users' feedback, as in the personalized optimization framework~\cite{simonetto2021personalized}.
	We also note that the unavailability of $(Q, R)$ renders methods relying on a learning phase for $(A, B)$ only (see, e.g.,~\cite{sforni2023on,sforni2024stability,song2024role}) inapplicable. %
	\oprocend
\end{remark} 
Our idea is to mimic~\eqref{eq:desired_update} by elaborating these finite-horizon approximations $\Jt(\g)$ according to an extremum-seeking perspective to compensate for the lack of knowledge about the gradient $\G(\g)$.

\section{\algo/: Algorithm Description and Convergence Properties}
\label{sec:algo}

\begin{figure}[htpb]
	\centering
	\includegraphics[scale=.9]{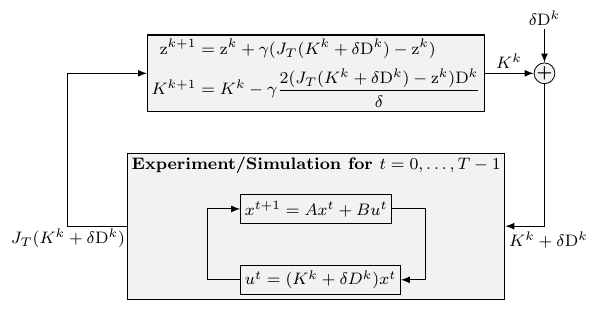}
	\caption{Block diagram representation of Algorithm~\ref{alg:algorithm}.}
	\label{fig:scheme}
\end{figure}

In this section, we present \algoext/ (\algo/), i.e., the novel data-driven method resumed in Algorithm~\ref{alg:algorithm} to iteratively address problem~\eqref{eq:problem} without the knowledge of the system and cost matrices $(A,B)$ and $(Q,R)$.
\begin{algorithm}[H]
	\begin{algorithmic}%
		\State \textbf{Initialization:} 
					$\f\ud 0 \in \R$, $\g\ud{0} \in \cK$. 
		\For{$\iter = 0, 1, 2 \ldots$}

		\State \textbf{Experiment phase (policy evaluation)}
		\State \hspace{2ex} Set the controller $\utime = (\gt + \delta\drt)\xtime$
		\State \hspace{2ex} Test $\xtimep = A\xtime + B\utime$ for $\timeid = 0,\dots,\Tsim -1$
		\State \hspace{2ex} Retrieve $\Jt(\gt+\delta\drt)$
		\State \textbf{Optimization phase (policy improvement)}
		\begin{subequations}\label{eq:system}
			\begin{align}
				\ftp &= \ft + \ssz(\Jt(\gt + \delta\drt) - \ft)
				\label{eq:system_ft}
				\\
				\gtp &= \gt - \ssz\frac{2(\Jt(\gt + \delta\drt) - \ft)\drt}{\delta}
				\label{eq:system_gt}
			\end{align}
		\end{subequations}
		\EndFor
		\caption{\algo/}
		\label{alg:algorithm}
	\end{algorithmic}
\end{algorithm}
Our algorithmic idea is to mimic the (model-based) gradient descent update~\eqref{eq:desired_update} through an extremum-seeking scheme.
To this end, at each iteration $\iter$, we perturb a given policy
gain $\gt$ obtaining $\gt + \delta\drt$, 
where $\delta > 0$ is an amplitude parameter and $\drt \in \R^{m \times n}$ is the so-called dither matrix.
The element $D_{ij}\ud\iter$ of $\drt$ is generated according to the sinusoidal law
\begin{align*}
	D_{ij}\ud\iter := \sin\left(\frac{2\pi\iter}{\Tijper} + \phi_{ij}\right),
\end{align*}
where $\Tijper \in \mathbb{Q}$ and $\phi_{ij} \in \mathbb{R}$ are the
period and the phase of component $(i,j)$, respectively, for all $(i,j) \in \until{\dimx} \times \until{\dimu}$.
Such a perturbed policy is used to implement the feedback control law $\utime = (\gt+ \delta\drt)\xtime$ and retrieve the corresponding finite-horizon cost $\Jt(\gt + \delta\drt)$ providing an approximation of the infinite-horizon one $\J(\g + \delta\drt)$.
This scenario may occur, for example, when a simulator of a complex system is available, but the analytical knowledge of the dynamics being implemented for the simulations is unavailable.
Hence, the finite-time truncation turns out to be crucial in avoiding experiments over infinite time horizons.
With $\Jt(\gt + \delta\drt)$ at hand, we perform the algorithm iteration detailed in~\eqref{eq:system}.
Specifically, the variable $\ft \in \R$ filters the variation of
$\Jt(\gt+\delta\drt)$ (see its update~\eqref{eq:system_ft}), while the
evolution of the gain matrix $\gt$ follows the extremum-seeking
update~\eqref{eq:system_gt}.
Unlike most existing algorithms (see, e.g.,~\cite{fazel2018global}) that rely on random perturbations, our extremum-seeking-based scheme benefits from the use of a deterministic dither matrix $\drt$, whose design conditions will be detailed in the following.
The combination of the dither matrix, the low-pass filter $\ft$, and ad-hoc deterministic experiments (see~\eqref{eq:exact_cost} and~\eqref{eq:cost_routine}) constitute the main novelties and the key advantages of the proposed scheme with respect to the existing methods.
A block diagram representation that graphically describes \algo/ is provided in Fig.~\ref{fig:scheme}.
Before establishing the convergence properties of \algo/, we need to ensure that the dither matrix is generated by following the orthonormality conditions detailed in the next assumption.
\begin{assumption}[Dither Frequencies Orthonormality]\label{ass:dither}
	The periods $\Toper, \dots, \Tnmper$ admit a least common multiple $\Tper \in \N$.
	Moreover, it holds
	\begin{subequations}\label{eq:frequencies}
		\begin{align}
			&
			\hspace{-.1cm}
			\sum_{\iter=1}^{\Tper} \sin\left(\tfrac{2\pi \iter}{ \Tpper} + \phi_{p}\right) \! = \! 0
			\label{eq:frequencies1}
			\\
			&
			\hspace{-.1cm}
			\sum_{\iter=1}^{\Tper} \sin\left(\tfrac{2\pi \iter}{ \Tpper} + \phi_{p}\right) \sin\left(\tfrac{2\pi \iter}{\Tpqer} + \phi_{q}\right) \! = \! \dfrac{\Tper}{2}
			\label{eq:frequencies2}
			\\
			&\hspace{-.1cm}
			\sum_{\iter=1}^{\Tper} \sin\!\left(\tfrac{2\pi \iter}{ \Tpper} \! + \! \phi_{p}\right)\! \sin\!\left(\tfrac{2\pi \iter}{\Tpqer}\! + \! \phi_{p}\right)\!\sin\!\left(\tfrac{2\pi \iter}{ \Trper} \! + \! \phi_{r}\right) \!=\! 0,\!\!\label{eq:frequencies3}
		\end{align}
	\end{subequations}
	for all $p,q,r \in \until{\dimu} \times \until{\dimx}$ such that $p \ne q$, $q \ne r$, and $p \ne r$. \oprocend
\end{assumption}
It is worth noting that a dither satisfying~\eqref{eq:frequencies1} and~\eqref{eq:frequencies2} can be interpreted as containing a set of orthogonal functions, in line with classic results from system identification theory.
Moreover, we also want~\eqref{eq:frequencies3} to be verified to improve the accuracy of the gradient estimation provided by the extremum seeking machinery. 
Conditions~\eqref{eq:frequencies1}-\eqref{eq:frequencies3}, which are in discrete time, are equivalent to those described in \cite[\S 2.1]{ariyur2003real} in continuous time.
Now, we are in the position to provide the main result of the paper, i.e., the convergence properties of \algo/.
\begin{theorem}[Convergence Properties of \algo/]\label{th:convergence}
	Consider \algo/ and let Assumptions~\ref{ass:A_B_Q_R},~\ref{ass:unknown}, and~\ref{ass:dither} hold.
	Then, for all $r > 0$ and $(\f\ud0,\g\ud0) \in \R \times \cK$, there exist $\bar{\ssz}, \bar{\delta}, a_0 > 0$, $a \in (0,1)$, and $\bar{\Tsim} \in \N$, such that, for all $\ssz \in (0,\bar{\ssz})$, $\delta \in (0,\bar{\delta})$, $\Tsim \ge \bar{\Tsim}$, the trajectories
	of~\eqref{eq:system} are bounded and satisfy
	\begin{subequations}\label{eq:thm_results}
		\begin{align}
			\gt + \delta\drt &\in \cK
			\label{eq:K_stabilizing}
			\\
			\norm{\gt - \gstar} &\leq a_0(1-\ssz a)^\iter + r,
			\label{eq:main_result}
		\end{align}
	\end{subequations}
	for all $\iter \in \N$.
	\oprocend
\end{theorem}
The proof of Theorem~\ref{th:convergence} is provided in Section~\ref{sec:proof}.
More in detail, the proof is based on the exploitation of Lyapunov stability and \emph{averaging theory} tools to prove that $(\f^\star,\gstar)$ is a semi-global practical exponentially stable equilibrium point of system~\eqref{eq:system} restricted to $\R \times \cK$, for a suitable $\f^\star \in \R$.

Theorem~\ref{th:convergence} requires an initial gain $\g\ud0$ that stabilizes the unknown pair $(A,B)$. 
We note that this does not necessarily require knowledge of $(A,B)$, see, e.g., the data-driven approach in~\cite{de2019formulas}, which computes stabilizing gains from system samples.  
\begin{remark}\label{rem:tuning}
	Theorem~\ref{th:convergence} formally guarantees that, for sufficiently small values of $\ssz$ and $\delta$, and sufficiently large values of $\Tsim$, \algo/ achieves the desired convergence properties.
	The proofs provide all the necessary steps to determine such values.
	However, in practice, obtaining them analytically may be complex. %
	As a result, one may prefer to determine $\bar{\ssz}$, $\bar{\delta}$, and $\bar{\Tsim}$ empirically, starting from reasonable initial guesses and then iteratively decreasing $\ssz$ and $\delta$, and increasing $\Tsim$, until the desired convergence properties are achieved. 
	\oprocend
\end{remark}

\section{Stability Analysis of \algo/}
\label{sec:analysis}
	
In this section, we perform the stability analysis of system~\eqref{eq:system} to prove Theorem~\ref{th:convergence}.
First, in Section~\ref{sec:FiniteSamples}, we perform a preliminary phase due to evaluate the approximation of the infinite-horizon gradient $\G(\g)$ using the finite-horizon cost $\Jt(\g)$.
In Section~\ref{sec:averaged_system}, by resorting to these approximations and an approach based on
\emph{averaging theory}, we characterize the stability and convergence properties of
the so-called \emph{averaged system} associated
to~\eqref{eq:system}.
With these results at hand, in Section~\ref{sec:proof}, we come back to the original time-varying system~\eqref{eq:system} and provide the proof of Theorem~\ref{th:convergence}.
Assumptions~\ref{ass:A_B_Q_R},~\ref{ass:unknown}, and~\ref{ass:dither} are valid throughout the entire section.

\subsection{Preliminary Approximation Results}
\label{sec:FiniteSamples}

Here, we provide two approximation results that will be used in the remainder of the analysis of system~\eqref{eq:system}.
First, we evaluate the approximation error due to using the truncated cost $\Jt(\g \!+\! \delta\drt)$ instead of the infinite-horizon one $\J(\g \!+\! \delta\drt)$.
\begin{lemma}[Truncated Cost Approximation Error]\label{lemma:finite_time_sim}
	For any $\acc > 0$ and compact set $\cS \subset \cK$, there exists $\bar{\Tsim} \in \N$ such that, for all $\Tsim \ge \bar{\Tsim}$, it holds  
	\begin{align}\label{eq:bound_sim}
		\|\J(\g) - \Jt(\g)\| \leq \acc,
	\end{align}
	for all $\g \in \cS$.\oprocend
\end{lemma}
The proof of Lemma~\ref{lemma:finite_time_sim} is provided in Appendix~\ref{sec:finite_time_proof}.

Second, we establish the gradient approximation properties obtained using the infinite-horizon cost samples $\J(\g+\dr\ud{1}),\dots,\J(\g+\dr\ud{\Tper})$ for any fixed (and stabilizing) gain $\g \in \cK$.
\begin{lemma}[Gradient Approximation Error]
	\label{lemma:estimation}
	For any compact set $\cS \subset \cK$, there exist $\err:\R^{\dimu \times \dimx}\to\R^{\dimu \times \dimx}$ and $\lippl > 0$ such that 
	\begin{subequations}
		\begin{align}
			\dfrac{2}{\delta \Tper}\!\!\sum_{\initer=\iter+1}^{\iter + \Tper}\!\!\! \J(\g + \delta \dr\ud\initer)\dr\ud\initer &= \G(\g) + \delta^2 \err(\g),
			\label{eq:result_gradient_estimation}
			\\
			\norm{\err(\g)} &\leq \lippl,\label{eq:ell_bound}
			\end{align}
	\end{subequations}
	for all $\iter \in \N$, $\delta \in (0,1]$, and $\g$ such that $\g + \delta\dr\ud\initer \in \cS$ for all $\initer \in \until{\Tper}$.
	\oprocend
\end{lemma}
The proof of Lemma~\ref{lemma:estimation} is provided in Appendix~\ref{sec:estimation proof}.

With these results at hand, we are able to study the stability properties of system~\eqref{eq:system} through the averaging theory.

\subsection{Averaged System Analysis}
\label{sec:averaged_system}

\begin{figure}[htpb]
	\centering 
	\includegraphics[scale=1]{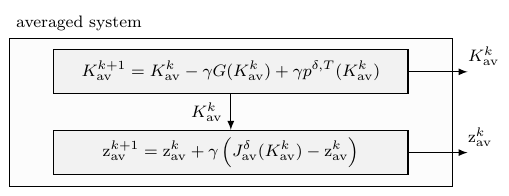}
	\caption{Block diagram representation of the averaged system~\eqref{eq:averaged_system}.}
	\label{fig:averaged}
\end{figure}

As shown in Section~\ref{sec:averaging_theory}, the averaged system associated to~\eqref{eq:system} is an
auxiliary dynamics derived by averaging the time-varying vector field
of~\eqref{eq:system} over time horizons of length equal to the period $\Tper$ (see Assumption~\ref{ass:dither}).
To properly write this system, given $\g \in \cK$ and $\z \in \R$, we consider the term $\sum_{\initer=\iter+1}^{\iter+\Tper}(\Jt(\g + \delta\dr\ud{\initer}) - \z)$ and add and subtract the infinite-horizon terms $\J(\g+\delta\dr\ud{\initer})$ with $\initer = 1,\dots,\Tper$, thus obtaining 
\begin{align}
	&
	\sum_{\initer=\iter+1}^{\iter+\Tper}(\Jt(\g + \delta\dr\ud{\initer}) - \z)\dr\ud\initer
	\notag\\
	&
	=\sum_{\initer=\iter+1}^{\iter+\Tper}(\J(\g + \delta\dr\ud{\initer}) - \z)\dr\ud{\initer} 
	\notag\\
	&\hspace{.4cm}
	+
	\sum_{\initer=\iter+1}^{\iter+\Tper}(\Jt(\g + \delta\dr\ud{\initer}) - \J(\g + \delta\dr\ud{\initer}))\dr\ud{\initer}
	\notag\\
	&\stackrel{(a)}{=} 
	\sum_{\initer=\iter+1}^{\iter+\Tper}\J(\g + \delta\dr\ud{\initer})\dr\ud{\initer} 
	\notag\\
	&\hspace{.4cm}
	+
	\sum_{\initer=\iter+1}^{\iter+\Tper}(\Jt(\g + \delta\dr\ud{\initer}) - \J(\g + \delta\dr\ud{\initer}))\dr\ud{\initer},\label{eq:sum_J}
\end{align}
where in $(a)$ we used the frequencies' property~\eqref{eq:frequencies1} to simplify the expression.
Hence, by applying Lemma~\ref{lemma:finite_time_sim}, Lemma~\ref{lemma:estimation}, and~\eqref{eq:sum_J}, the averaged system
associated to~\eqref{eq:system} reads as%
\begin{subequations}\label{eq:averaged_system}
	\begin{align}
		\fatp &= \fat + \ssz\left(\Ja(\gat) - \fat\right)\label{eq:averaged_system_ft}
		\\
		\gatp &= \gat - \ssz\G(\gat) + \ssz\app(\gat),\label{eq:averaged_system_gt}
	\end{align}
\end{subequations}
where $\Ja: \cK \to \R$ and $\app: \cK \to \cK$ are defined as 
\begin{subequations}
	\begin{align}
		\Ja(\g) &:= \frac{1}{\Tper}\sum_{\initer = \iter+1}^{\iter+\Tper}\Jt(\g+\delta\dr\ud\initer)
		\\
		\app(\g) &:=-\dfrac{2}{\delta \Tper}\!\!\sum_{\initer=\iter+1}^{\iter + \Tper}\!\!\! \left(\J_\Tsim(\g \! + \! \delta \dr\ud\initer) \! - \! \J(\g \! + \! \delta \dr\ud\initer)\right)\dr\ud\initer\!
		\notag\\
		&\hspace{.4cm}
		-\delta^2 \err(\g).
		\label{eq:app}
	\end{align}
\end{subequations}
As graphically highlighted in Fig.~\ref{fig:averaged}, we remark that the averaged scheme~\eqref{eq:averaged_system} is a cascade system.
The next lemma provides the convergence properties of the averaged system~\eqref{eq:averaged_system}.
To this end, we introduce the candidate Lyapunov function $V: \R \times \cK \to \R_+$ defined as
\begin{align}
	\V(\f,\g) &:=  \frac{1}{2 \mav}\norm{\f}^2+ \J(\g) - \J(\gstar),\label{eq:V_definition}
\end{align}
where $\mav \ge 1$ will be fixed in the next lemma.
\begin{lemma}[Averaged System Stability]\label{lemma:averaged}
	Consider~\eqref{eq:averaged_system}.
	Then, for all $(\fa\ud0,\ga\ud0) \in \R \times \cK$ and $\rav > 0$, there exist $\bar{\ssz}_1,\bar{\delta}_1, a, \bar{\Tsim} > 0$ and $\bar{\mav} \ge 1$ such that, for all $\ssz \in (0,\bar{\ssz}_1)$, $\delta \in (0,\bar{\delta}_1)$, $\Tsim \ge \bar{\Tsim}$, and $\mav \ge \bar{\mav}$, it holds 
	\begin{align}
		&\V(\fat - \Ja(\gat),\gat) 
		\notag\\
		&\leq (1 - \ssz a)^\iter \V(\fa\ud0 - \Ja(\ga\ud0),\ga\ud0) + \rav,\label{eq:exp_bound}
	\end{align}
	for all $\iter \in \N$.
	\oprocend
\end{lemma}
The proof of Lemma~\ref{lemma:averaged} is provided in Appendix~\ref{sec:proof_lemma_averaged}.

\subsection{Proof of Theorem~\ref{th:convergence}}
\label{sec:proof}

The proof relies on the application of Theorem~\ref{th:averaging}
(cf. Section~\ref{sec:averaging_theory}) to system \eqref{eq:system}.
Then, in order to apply Theorem~\ref{th:averaging}, we need to (i) choose the design parameters $c_1, \rho > 0$ bounding the initial and final values of $\V$, respectively, and (ii) satisfy the conditions required by Assumptions~\ref{hyp:fa},~\ref{hyp:bounds}, and~\ref{hyp:StabilityAv}. 
By~\cite[Lemma~3.8]{bu2020lqr}, we recall that there exists $\gl > 0$ such that 
\begin{align}
	\gl\norm{\g - \gstar}^2 \leq \J(\g) - J(\gstar),\label{eq:lower_bound}
\end{align}
for all $\g \in \cK$.
Therefore, by looking at the statement of Theorem~\ref{th:averaging} and given the desired final radius $r$, we set $\rho \in (0,\sqrt{r\gl}]$.
In order to set the initial radius, we need to find a bound for $\delta$ such that $\g\ud0 + \delta\drt$ is stabilizing for all $\iter \in \N$.
To this end, we note that $\g\ud0 \in \cK$, $\cK$ is open~\cite[Lemma~IV.3]{bu2020topological}, and $\drt$ is bounded for all $\iter \in \N$.
Hence, there exists $\bar{\delta}_0 > 0$ such that $\g\ud0 + \delta\drt \in \cK$ for all $\delta \in [0,\bar{\delta}_0]$ and $\iter \in \until{\Tper}$.
Now, we arbitrarily choose $c_1 > \max_{\delta \in [0,\bar{\delta}_0]}V_1(\f\ud0 - \Ja(\g\ud0),\g\ud0)$ and, thus, we note that $c_1 \ge \max_{\delta \in [0,\bar{\delta}]}\V(\f\ud0 - \Ja(\g\ud0),\g\ud0)$ for all $\mav \ge 1$ (see the definition of $\V$ in~\eqref{eq:V_definition}).
Once the initial and final radius $c_1$ and $\rho$ have been chosen, let us check Assumptions~\ref{hyp:fa},~\ref{hyp:bounds}, and~\ref{hyp:StabilityAv}. %
First, Assumption~\ref{hyp:fa} is trivially satisfied because the dither signals are $\Tper$-periodic (cf. Assumption~\ref{ass:dither}).  
Second, 
we remark that~\eqref{eq:system} and its corresponding averaged system~\eqref{eq:averaged_system} are continuous over the set $\{(\f,\g)\in \R \times \cK \mid \g + \delta\drt \in \cK \text{ for all } \iter\in\N\}$, as required by Assumption~\ref{hyp:bounds}.
For this reason, let us choose $\delta$ such that the level set $\levset_{c_1} := \{(\f,\g) \in \R \times \cK \mid V_1(\z- \Ja(\g),\g) \leq c_1\}$ of $V_1$ (i.e, the function $\V$ with $\mav =1$, see~\eqref{eq:V_definition}) is contained into $\{(\f,\g)\in \R \times \cK \mid \g + \delta\drt \in \cK \text{ for all } \iter\in\N\}$.
To this end, by looking at the definition of $\V$ (cf.~\eqref{eq:V_definition}), we note that
\begin{align*}
	V_1(\f - \Ja(\g),\g) \leq c_1 \implies \J(\g) - \J(\gstar) \leq c_1,
\end{align*}
independently on the choice of $\delta$.
In turn, the result~\eqref{eq:implication} implies
\begin{align} 
	(\f,\g) \in \levset_{c_1} \implies \g \in \cK.\label{eq:implication}
\end{align}
Moreover, we recall that $\cK$ is open~\cite[Lemma~IV.3]{bu2020topological} and $\drt$ is bounded for all $\iter\in\N$.
Then, we guarantee the existence of $\bar{\delta}_2 > 0$ such that, for all $\delta \in (0,\min\{\bar{\delta}_0,\bar{\delta}_2\})$, it holds $\g + \delta\drt \in \cK$ for all $\g$ satisfying $\J(\g) - \J(\gstar) \leq c_1$ and $\iter\in\N$.
With these results at hand, Lemma~\ref{lemma:averaged} ensures the existence of $\bar{\ssz}_1, \bar{\delta}_1 > 0$, $\bar{\mav} \ge 1$, and $\bar{\Tsim} \in \N$ such that, by setting $\ssz \in (0,\bar{\ssz}_1)$, $\delta \in (0,\bar{\delta})$ with $\bar{\delta} := \min\{\bar{\delta}_0,\bar{\delta}_1,\bar{\delta}_2\}$, $\Tsim \ge \bar{\Tsim}$ and $\mav \ge \bar{\mav}$, $\V$ achieves the convergence properties~\eqref{eq:exp_bound} along the trajectories of the averaged system~\eqref{eq:averaged_system} and, thus, Assumption~\ref{hyp:StabilityAv} is satisfied.
Hence, we are entitled to apply Theorem~\ref{th:averaging} which, for all $(\f\ud0,\g\ud0) \in \Omega_{c_0}$, ensures the existence of $\bar{\ssz} > 0$ such that, for all $\ssz \in (0,\bar{\ssz})$, the trajectories of~\eqref{eq:system} satisfy
\begin{subequations}
	\begin{align}
		\V(\ft   -  \J(\gt),\gt)  &\leq  c_1
		\label{eq:lypaunov_bounded} 
		\\
		\V(\ft  -  \J(\gt),\gt)  &\leq  (1  -  \ssz a)^\iter \V(\f\ud0  -  \J(\g\ud0),\g\ud0)
		 +  \rho,
		\label{eq:final_cost}
	\end{align}
\end{subequations}
for all $\iter \in \N$.
The proof of~\eqref{eq:K_stabilizing} follows by combining~\eqref{eq:implication} and~\eqref{eq:lypaunov_bounded},
while the proof of~\eqref{eq:main_result} follows by
combining~\eqref{eq:lower_bound}, $\J(\g) - \J(\gstar) \leq \V(\f -
\J(\g),\g)$,~\eqref{eq:final_cost}, the
choice of $\rho$, and by setting $a_0 := \sqrt{\V(\f\ud0  -  \J(\g\ud0),\g\ud0)/\gl}$.
\section{Numerical Simulations: Control of a Doubly Fed Induction Motor}
\label{sec:numerical_simulations}

\begin{figure}[H]
	\centering
	\includegraphics[scale=1]{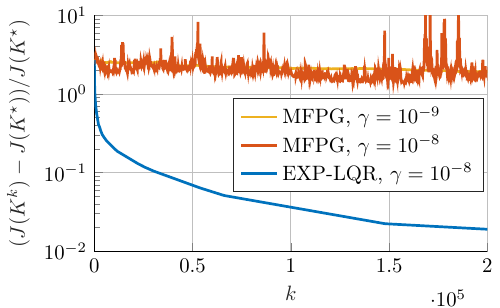}
	\caption{Comparison between \algo/ and MFPG~\cite{fazel2018global} in terms of the relative cost error $(\J(\gt) - \J(\gstar))/\J(\gstar)$ along the algorithms' iteration $\iter$.}
	\label{fig:cost_error}
\end{figure}

\begin{figure}[htpb]
	\centering
	\includegraphics[scale=1]{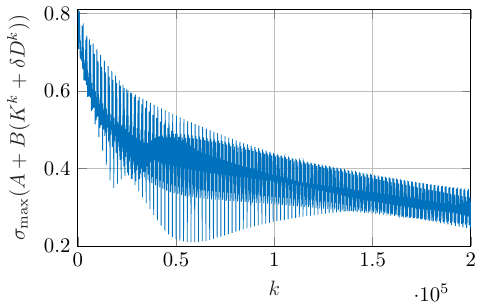}
	\caption{Evolution of the maximum (in absolute value) closed-loop matrix eigenvalue $\sigma_{\text{max}}(A + B(\gt + \delta\drt))$ along the algorithm iteration $\iter$.}
	\label{fig:max_eig}
\end{figure}

\begin{figure}[htpb]
	\centering
	\subfloat[][Iteration $\iter = \itercodea$.]
	   {\makebox[0.45\linewidth]{\includegraphics[scale=.4999]{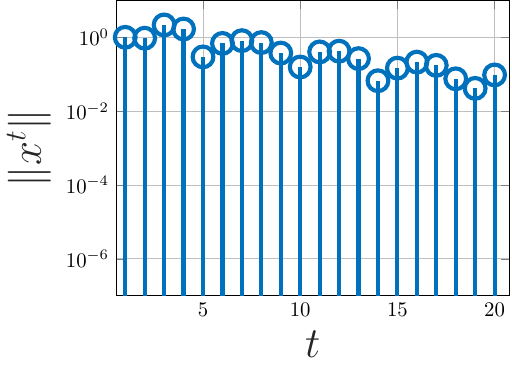}}}
	   \quad\quad
	\subfloat[][Iteration $\iter = \itercodeb$.]
	   {\makebox[0.45\linewidth]{\includegraphics[scale=.4999]{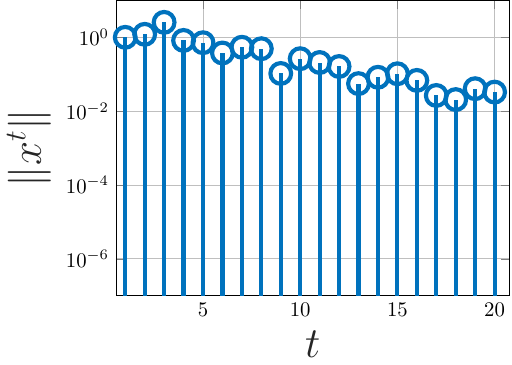}}}
	\\
	\subfloat[][Iteration $\iter = \itercodec$.]
	   {\makebox[0.45\linewidth]{\includegraphics[scale=.4999]{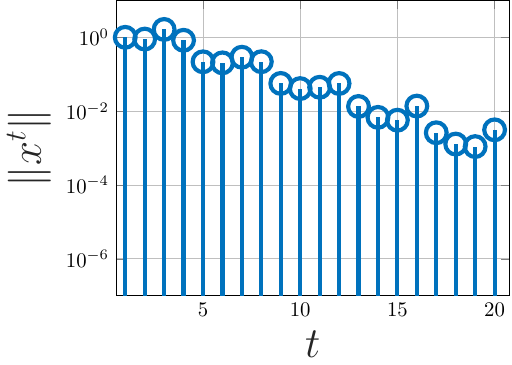}}}
	   \quad\quad
	\subfloat[][Iteration $\iter = \itercoded$.]
	   {\makebox[0.45\linewidth]{\includegraphics[scale=.4999]{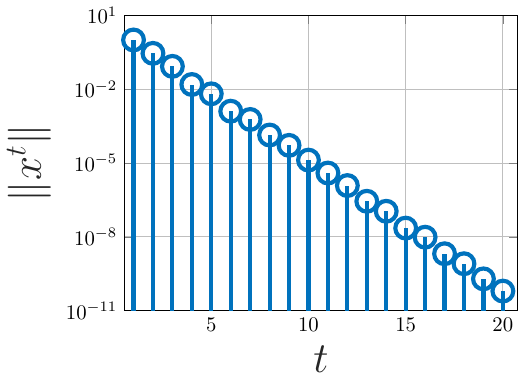}}}
	\caption{System evolutions in the interval $[0,\Tsim]$ along different algorithm iterations $\iter$.}
	\label{fig:traj}
\end{figure}

In this section, we numerically test the effectiveness of \algo/ and compare it with the Model-Free Policy Gradient (MFPG) method proposed in~\cite{fazel2018global}.
To this end, we consider a forward Euler discretization of the continuous-time linear model provided by~\cite{leonhard2001control} for a Doubly Fed Induction Motor (DFIM) operating at constant speed.
Namely, we consider the discrete-time linear system
\begin{align}\label{eq:system_sim}
	\xtimep = \underbrace{(I + \ts \Ac)}_{A}\xtime + \underbrace{\ts \Bc}_{B}\utime,
\end{align}
where $\ts = \tscode$ is the adopted sampling period, while $x, u \in \R^{4}$ are the state and input variables and are defined as
\begin{align*}
	x := \begin{bmatrix}
		i_{1u}& i_{1v}& i_{2u}& i_{2v}
	\end{bmatrix}\T, \hspace{.1cm} u := \begin{bmatrix}u_{1u}& u_{1v}& u_{2u}& u_{2v}\end{bmatrix}\T\!\!,
\end{align*}
where $i_{1u}, i_{1v} \in \R$ are the stator currents and $i_{2u}, i_{2v} \in \R$ are the rotor currents, while $u_{1u}, u_{1v} \in \R$ are the stator voltages and $u_{2u}, u_{2v} \in \R$ the rotor voltages. 
Finally, $\Ac \in \R^{4 \times 4}$ and $\Bc \in \R^{4\times 4}$ represent the state and input matrices of the continuous-time model, respectively, and are defined as
\begin{align*}%
	\Ac &:= \frac{1}{\bar{L}}
	\begin{bmatrix}
		-L_2R_1 &  -a + b & L_mR_2 & b_2\\
			a - b & -L_2R_1 & -b_2 & -L_mR_2 \\
		L_mR_1 & -b_1 & -L_1R_2 & -a-b_{12}\\
		b_1 & L_mR_1 &  a+b_{12}& -L_1R_2
	\end{bmatrix}\\
	\Bc &:= \frac{1}{\bar{L}}
	\begin{bmatrix}
		L_2 & 0 & -L_m & 0\\
		0 & L_2 & 0 & -L_m\\
		-L_m & 0 & L_1 & 0\\
		0 & -L_m & 0 & L_1
	\end{bmatrix},
\end{align*}
where 
\begin{align*}
	\bar{L}&\coloneqq L_1L_2-L_m^2
	,\quad
	a := \bar{L}\omega_0, \quad b\coloneqq L_m^2\omega_r\\
	b_{12} &:= L_1L_2\omega_r, \quad b_1 \coloneqq L_1L_m\omega_r,\quad b_2 \coloneqq L_2L_m\omega_r.
\end{align*}
More in detail, the parameters $R_1$ and $R_2$ correspond to the resistances of the stator and rotor, while the parameters $L_1$, $L_2$, and $L_m$ refer to the stator and rotor self-inductances, and the mutual inductance, respectively. 
Lastly, $\omega_r$ and $\omega_0$ denote the electrical angular velocities of the rotor and the rotating reference frame, respectively, which are assumed constant.
We adopt the physical parameters used by~\cite{borghesi2024mr} about the same model, and we report them in Table~\ref{table:params}.
These parameters make the discrete-time pair $(A,B)$ controllable as required by Assumption~\ref{ass:A_B_Q_R}.
\addtocounter{table}{0}
\begin{table}[htpb]
\caption{Physical parameters of the DFIM~\cite{borghesi2024mr}.}
\label{table:params}
\centering
\begin{tabular}{c | c | c | c}
	Parameter & Value & Parameter & Value \\
	\hline
	$L_1$ & $0.02645$ & $R_1$ [$\Omega$]& $0.036$ \\
	$L_2$ & $0.0264$& $R_2$ [$\Omega$]& $0.038$ \\
	$L_m$ & $0.0257$& $\omega_0$ [rad/s]& $2\pi70.8$\\
	$p$ & $3$ & $\omega_r$ [rad/s]& $2\pi62$  \\
\end{tabular}
\end{table}
For the cost matrices $Q \in \R^{4 \times 4}$ and $R \in \R^{4 \times 4}$, we randomly generate them to ensure they are symmetric, with eigenvalues lying within the interval $(0, \Qcode)$ thus satisfying Assumption~\ref{ass:A_B_Q_R}.
We empirically tune the algorithm parameters as $\ssz = \sszcode$, $\delta = \deltacode$, and $\Tsim = \Tsimcode$ (see Remark~\ref{rem:tuning}).
Namely, the algorithm needs to perform $4$ experiments or simulations given by $\Tsimcode$ samples per iteration to retrieve the truncated cost $\Jt(\gt+\delta\drt)$ (see~\eqref{eq:exact_cost} and~\eqref{eq:cost_routine}).
As for the generation of the dither matrix $\drt \in \R^{4 \times 4}$, we ordered the pairs $(i,j) \in \until{4} \times \until{4}$ with indices $p=1, \dots, 16$ and chosen $\mathbf{\iter}_{p,\text{prd}} = \Tpercode \times 2^{(1-p)/2}$ and $\phi_p = 0$ for $p$ odd, while $\mathbf{\iter}_{p,\text{prd}} = \mathbf{\iter}_{p-1,\text{prd}}$ and $\phi_p = \pi/2$ for $p$ even.
This choice ensures that Assumption~\ref{ass:dither} is satisfied with period $\Tper = \Tpercode$ by directly applying~\cite[Theorem]{Knapp2009Sines}.
As for MFPG, we adopt a tuning analogous to that of \algo/.
Specifically, we run MFPG using the same step size, perturbation amplitude, number of samples per experiment, number of experiments per iteration, and variables' initialization used in \algo/.
For the sake of fairness, we also report the results achieved by MFPG with $\gamma = 10^{-8}$.
Fig.~\ref{fig:cost_error} shows the evolution of the relative cost error $(\J(\gt) - \J(\gstar))/\J(\gstar)$ along the algorithms' iterations $\iter$ in logarithmic scale.
As predicted by Theorem~\ref{th:convergence}, Fig.~\ref{fig:cost_error} 
shows that \algo/ asymptotically converges in a neighborhood of the optimal gain $\gstar$. 
Moreover, Fig.~\ref{fig:cost_error} shows that \algo/ outperforms MFPG in terms of both convergence speed and final accuracy.	
In Fig.~\ref{fig:max_eig}, we show the evolution of $\sigma_{\text{max}}(A + B(\gt + \delta\drt))$ along the algorithm iterations $\iter$, where, given a generic square matrix $M \in \R^{n \times n}$, the symbol $\sigma_{\text{max}}(M)$ denotes the maximum (in absolute value) eigenvalue of $M$.
In particular, Fig.~\ref{fig:max_eig} shows that $\sigma_{\text{max}}(A + B(\gt + \delta\drt))$ never reaches the unitary value.
Namely, as predicted by Theorem~\ref{th:convergence} (cf.~\eqref{eq:K_stabilizing}), we always test the system through a stabilizing state-feedback controller $\utime = (\gt + \delta\drt)\xtime$.
Finally, in Fig.~\ref{fig:traj}, we show the evolution of the norm of the state trajectory $\norm{\xtime}$ of system~\eqref{eq:system_sim} in four simulations (each composed of $\Tsim = \Tsimcode$ samples) performed at different algorithm iterations $\iter$ to retrieve the truncated cost $\Jt(\gt+\delta\drt)$.
In particular, Fig.~\ref{fig:traj} shows that the trajectories of system~\eqref{eq:system_sim} (controlled with $\utime = (\gt + \delta\drt)\xtime$) exponentially converge to the origin quicker and quicker as the iteration index $\iter$ increases since we are iteratively reducing the absolute values of the eigenvalues of the gain closed-loop matrix $(A + B\gt)$ (see also Fig.~\ref{fig:max_eig}). 

\section{Conclusions}

We proposed \algo/, i.e., a novel data-driven method able to iteratively find the state feedback gain matrix solving a Linear Quadratic Regulator problem.
\algo/ does not need the direct knowledge of the system and cost matrices.
Indeed, given an oracle able to provide a finite-time truncation of the LQR cost, our method refines its estimate according to a mechanism based on extremum-seeking.
We analyzed the resulting time-varying algorithm by exploiting system theory tools based on Lyapunov stability and averaging theory.
Specifically, we guaranteed that our algorithm exponentially converges to an arbitrarily small ball containing the optimal gain matrix.
We tested the proposed solution with numerical simulations involving the control of an inductance motor.

\appendix

\subsection{Proof of Theorem~\ref{th:averaging}} 
\label{sec:proof_averaging_generic}

Since Assumption~\ref{hyp:StabilityAv} characterizes the evolution of $V$ along the trajectories $\{\state\iter\av\}_{\iter\in \N}$ of the averaged system~\eqref{eq:averaged_system}, the idea of the proof is to bound the distance $\|\state\iter - \state\iter\av\|$ to characterize the evolution of $V$ along the trajectories $\{\state\iter\}_{\iter\in \N}$ of the original time-varying system~\eqref{eq:plant}.
To this end, we introduce $\upsilon: \R^n \times \N \to \R^{n}$ defined as
\begin{equation}
\label{eq:upsilon}
    \upsilon(\state{}\av,\iter) := \sum_{\initer=0}^{\iter-1}\left(f(\state{}\av,\initer)-f\av(\state{}\av)\right).
\end{equation}
By using~\eqref{eq:fa} and~\eqref{eq:upsilon}, the evolution of $\upsilon$ reads as
\begin{align}
&\upsilon(\state\iterp\av,\iterp)- \upsilon(\state\iter\av,\iter) 
\notag\\
&= 
f(\state\iterp\av,\iter)-f\av(\state\iterp\av)+ \upsilon(\state\iterp\av,\iter) - \upsilon(\state\iter\av,\iter).\label{eq:upsilon_update}
\end{align}
Let us recall that $V(\stateinit) \leq c_0$ and that $\rho < c_0 < c_1$ by assumption.
Then, let us arbitrarily choose $\epsilon \in (0,\min\{c_1 - c_0,\rho\})$, and $\rav \in (0,\min\{\rho - \epsilon,c_0\})$.
As it will become clearer later, $\epsilon$ represents the maximum difference between $V(\state{\iter})$ and $V(\state{\iter}\av)$, where $\rav$ defines the level set of $V$ in which we enforce the convergence of the averaged state $\state\iter\av$.
(cf. Assumption~\ref{hyp:StabilityAv}).
Under the assumption of $\state{\iter} \in \Omega_{c_1}$ for all $\iter\in\N$ (later verified by a proper selection of $\ssz$), we use the compactness of the set $\Omega_{c_1}$ (cf. Assumption~\ref{hyp:StabilityAv}) and the continuity properties over $\Omega_{c_1} \subseteq \cX$ (cf. Assumption~\ref{hyp:bounds}) to ensure the existence of $L > 0$ such that  
\begin{subequations}\label{eq:bounds}
	\begin{align}
		\norm{f(\state{},\iter)} &\leq L
		, \quad 
		\norm{f\av(\state{})} \leq L 
		\\
		\norm{\dfrac{\partial f(\state{},\iter)}{\partial \state{}}} &\leq L
		,\quad
		\norm{\nabla f\av(\state{})} \leq L 
		\\
		\norm{\nabla V(\state{})} &\leq L,
	\end{align}
\end{subequations}
for all $\state{} \in \Omega_{c_1}$ and $\iter\in\N$.
In turn, the bounds~\eqref{eq:bounds} lead to
\begin{subequations}\label{eq:bounds_phi}
	\begin{align}
		\norm{\upsilon(\state{},\iter)} &\leq 2 L \Tper\label{eq:bound_v}
		\\
		\norm{f(\state{},\iter)-f(\state\prime,\iter)} &\leq  L \norm{\state{}-\state\prime}
		\\
		\norm{f\av(\state{})-f\av(\state\prime)} &\leq  L \norm{\state{}-\state\prime}
		\\
		\norm{\upsilon(\state{},\iter) - \upsilon(\state\prime,\iter)} &\leq 2L \Tper  \norm{\state{}-\state\prime}
		\\
		\norm{f\av(\state{})} &\leq L
		\\
		V(\state{}) - V(\state\prime) &\leq L\norm{\state{}-\state\prime},\label{eq:lip_lyap}
	\end{align}
\end{subequations}
for all $\state{}, \state\prime \in \Omega_{c_1}$ and $\iter\in\N$.
Now, let us introduce $\zeta^\iter \in \R^{n}$ defined as
\begin{align}
	\zeta^\iter := \state\iter\av+\ssz \upsilon(\state\iter\av,\iter).\label{eq:zeta}
\end{align} %
By algebraically rearranging the terms, we can write 
\begin{align*}
	\state\iter-\zeta^\iter &=  \sum_{\initer = 0}^{\iter-1} \left((\state{\initer+1}-\state{\initer})-(\zeta^{\initer+1}-\zeta^\initer)\right).
\end{align*}
Now let us add $\pm\ssz \sum_{\initer = 0}^{\iter-1}(f(\zeta\ud\initer,\initer)  + f(\state{\initer}\av,\initer))$ in the above equation and use~\eqref{eq:upsilon_update} to get
\begin{align}
	\state\iter-\zeta^\iter 
	&=  \ssz \sum_{\initer = 0}^{\iter-1}(f(\state{\initer},\initer)-f(\zeta^\initer,\initer))
	\notag\\
	&\hspace{.4cm}
	+\ssz\sum_{\initer = 0}^{\iter-1}(f(\zeta^\initer,\initer)-f(\state{\initer}\av,\initer))
	\notag\\
	&\hspace{.4cm}
	-\ssz\sum_{\initer = 0}^{\iter-1}(f(\state{\initer+1}\av,\initer)-f(\state{\initer}\av,\initer))
	\notag\\
	&\hspace{.4cm}
	+\ssz\sum_{\initer = 0}^{\iter-1}(f\av(\state{\initer+1}\av)-f\av(\state{\initer}\av))
	\notag\\
	&\hspace{.4cm}
	-\ssz\sum_{\initer = 0}^{\iter-1} ( \upsilon(\state{\initer+1}\av,\initer) - \upsilon(\state{\initer}\av,\initer)).\label{eq:chi_minus_zeta}
\end{align}
By combining~\eqref{eq:chi_minus_zeta} with~\eqref{eq:plant},~\eqref{eq:average}, and~\eqref{eq:bounds_phi}, we can write
\begin{align}
	\norm{\state\iter-\zeta^\iter} \!\le\!  \ssz L \sum_{\initer = 0}^{\iter-1}\norm{\state{\iter}-\zeta^\initer}+\ssz^2 L^2 2 \left( 1 + 2  \Tper  \right) \iter.\!\label{eq:chi_zeta_bound}
\end{align}
Note that 
\begin{align}
	\sum_{\initer = 0}^{\iter-1}  \ssz L k  \exp\left(- \ssz L k\right) \le \sum_{\initer = 0}^{\infty}  \ssz L \initer \exp\left(- \ssz L \initer \right) = 1.\label{eq:exp_one}
\end{align}
By combining~\eqref{eq:exp_one} and the discrete Gronwall inequality (see~\cite{Popenda1983Discrete,holte2009discrete}), we are able to bound~\eqref{eq:chi_zeta_bound} as
\begin{align*}
	\norm{\state\iter-\zeta^\iter}
	\!\le \! \ssz^2 L^2 2 \left( 1 \! + \! 2  \Tper  \right) \iter \! + \!  
	\ssz L 2 \left( 1 \! + \! 2  \Tper  \right)  \exp\left(\ssz L \iter \right)\!.
\end{align*}
By combining the latter with the definition of $\zeta$ (cf.~\eqref{eq:zeta}) and the triangle inequality, we get
\begin{align}
	\norm{\state\iter-\state\iter\av} 
	&\leq \ssz^2 L^2 2 \left( 1 + 2  \Tper  \right) \iter 
	\notag\\
	&\hspace{.4cm}
	+ \ssz L 2 \left( 1 + 2  \Tper  \right)  \exp\left(\ssz L \iter \right) \! + \! \ssz\norm{\upsilon(\iter,\state\iter\av)}
	\notag\\
	&\stackrel{(a)}{\leq} \ssz^2 L^2 2 \left( 1 + 2  \Tper  \right) \iter
	\notag\\
	&\hspace{.4cm}
+\ssz L 2 \left( 1 + 2  \Tper  \right)  \exp\left(\ssz L \iter \right) \! + \! \ssz 2L \Tper.\!\!\label{eq:bound_x_xa}
\end{align}
where in $(a)$ we use~\eqref{eq:bound_v} to bound $\norm{\upsilon(\iter,\state\iter\av)}$.
Then, we set $\theta_{\rav} \in \R$ such that
\begin{align}
	\theta_{\rav} \ge -\frac{1}{a}\ln\left(\tfrac{\rav}{c_0}\right).
	\label{eq:theta_star}
\end{align}
Now, we want to impose the $\epsilon$-closeness between the trajectories of system~\eqref{eq:plant} and its averaged version~\eqref{eq:average}.
To this end, by looking at the bound in~\eqref{eq:bound_x_xa}, we introduce 
\begin{subequations}
	\begin{align}
		\bar{\ssz}_2 &:= \tfrac{\epsilon/(3L)}{L^2 2 \left( 1 + 2  \Tper  \right) \theta_{\rav}} 
		,\quad
		\bar{\ssz}_3 := \tfrac{\epsilon/(3L)}{2L (1+2\Tper )\exp(L \theta_{\rav})}
		\\
		\bar{\ssz}_4 &:= \tfrac{\epsilon/(3L)}{2 L \Tper }
		,\quad
		\hspace{1.3cm}\bar{\ssz} := \min\{\bar{\ssz}_1, \bar{\ssz}_2,\bar{\ssz}_3,\bar{\ssz}_4\}.\label{eq:bar_gamma_1}
	\end{align}
\end{subequations}
Subsequently, we pick $\ssz \in (0,\bar{\ssz})$ such that $\bar{\iter} := \frac{\theta_{\rav}}{\ssz}\in\N$. 
This can be done without loss of generality since $\theta_{\rav}$ is a design parameter.
Then, the definition of $\bar{\ssz}$ (cf.~\eqref{eq:bar_gamma_1}) and the inequality in~\eqref{eq:bound_x_xa} lead to the bound
\begin{align}
	\norm{\state\iter-\state\iter\av} \leq \epsilon/L,\label{eq:closeness}
\end{align}
for all $\iter \in \{0,\dots,\bar{\iter}\}$.
Then, for all $\iter \in \{0,\dots,\bar{\iter}\}$, we add and subtract $V(\state\iter\av)$ to $V(\state\iter)$ and write
\begin{align}
	V(\state\iter) &= V(\state\iter\av) + V(\state\iter) - V(\state\iter\av) 
	\notag\\
	&\stackrel{(a)}{\leq} 
	c_0 + L\norm{\state\iter - \state\iter\av}
	\notag\\
	&\stackrel{(b)}{\leq} 
	c_0 + \epsilon
	\stackrel{(c)}{\leq} 
	c_1,
	\label{eq:x_stay_c1}
\end{align}
where in $(a)$ we use the fact that $\state\iter\av \in \Omega_{c_0}$ for all $\iter\ge 0$ (see~\eqref{eq:StabilityAv} by Assumption~\ref{hyp:StabilityAv}) and the bound~\eqref{eq:lip_lyap}, in $(b)$ we use the bound~\eqref{eq:closeness}, while in $(c)$ we use the fact that $\epsilon \leq c_1 - c_0$.
Therefore, the bound~\eqref{eq:x_stay_c1} allows us to claim that ${\state{}}^\iter \in \Omega_{c_1}$ for all $\iter \in \{0,\dots,\bar{\iter}\}$, i.e., we have verified that the bounds~\eqref{eq:bounds_phi} can be used in the interval $\{0,\dots,\bar{\iter}\}$.
Further, the exponential law~\eqref{eq:StabilityAv} and the expression of $\theta_{\rav}$ (cf.~\eqref{eq:theta_star}) ensure that it holds
\begin{align}
	V(\state\iter\av) \leq \rav,\label{eq:exponential_tstar}
\end{align}
for all $\iter \ge \bar{\iter}$.
By adding $\pm V(\state{\bar{\iter}}\av)$ to $V(\state\iter)$, we get
\begin{align}
	V(\state{\bar{\iter}}) &= V(\state{\bar{\iter}}\av) + V(\state{\bar{\iter}}) - V(\state{\bar{\iter}}\av)
	\notag\\
	&
	\stackrel{(a)}{\le} \rav + L\norm{\state{\bar{\iter}} - \state{\bar{\iter}}\av} 
	\stackrel{(b)}{\le} 
	\rav + \epsilon
	\stackrel{(c)}{\le} \rho,\label{eq:ball}
\end{align}
where in $(a)$ we combined~\eqref{eq:lip_lyap} and~\eqref{eq:exponential_tstar}, in $(b)$ we used~\eqref{eq:closeness}, while $(c)$ uses the choice of $\rav \leq \rho - \epsilon$.
We remark that the inequality~\eqref{eq:ball} also guarantees that $\state{\bar{\iter}} \in \Omega_{c_0}$ since $\rho \leq c_0$.
Next, in order to show that $\state\iter \in \Omega_{\rho}$ for all $\iter \ge \bar{\iter}$, we divide the set of natural numbers in intervals as 
$\mathbb{N} = \{0,\dots,\bar{\iter}\}\cup \{\bar{\iter},\dots,2\bar{\iter}\} \cup \dots$.
Define $\psi\av(\state{\bar{\iter}},\iter+ \bar{\iter})$ as the solution to~\eqref{eq:average} for $\state{0}\av = \state{\bar{\iter}}$ and $\iter \in \{0,\dots,\bar{\iter}\}$.
Thus, at the beginning of the time interval $\{\bar{\iter},\dots,2\bar{\iter}\}$, the initial condition of the trajectory of~\eqref{eq:average} coincides with the one of $\psi\av(\state{\bar{\iter}},\iter+\bar{\iter})$ and lies into $\Omega_{\rho} \subseteq \Omega_{c_0}$.
Thus, we apply the same arguments above to guarantee that, for any $\ssz \in (0,\bar{\ssz})$, it holds
\begin{align*}
	\norm{\state{k+\bar{\iter}}-\psi\av(\state{\bar{\iter}},\iter+\bar{\iter})} \leq \epsilon/L,
\end{align*}
for all $\iter \in \{0,\dots,\bar{\iter}\}$.
Moreover, with the same arguments, it holds $\psi\av(\state{\bar{\iter}},2\bar{\iter}) \in \Omega_{\rav}$.
Then, in light of Assumption~\ref{hyp:StabilityAv}, we guarantee that the averaged system~\eqref{eq:average} cannot escape from the set $\Omega_{\rho}$, namely, for all $\iter \in \{\bar{\iter},\dots,2\bar{\iter}\}$, it holds
\begin{align*}
	\state\iter\av \in \Omega_{\rho}.
\end{align*}
Thus, we get $\state\iter \in \Omega_{\rho}$ for all $\iter \in \{\bar{\iter},\dots,2\bar{\iter}\}$.
The proof follows by recursively applying the above arguments for each time interval $\{j\bar{\iter},\dots,(j+1)\bar{\iter}\}$ with $j = 2, 3, \dots$.

\subsection{Proof of Lemma~\ref{lemma:finite_time_sim}}
\label{sec:finite_time_proof}

We observe that
\begin{align*}
	\lim_{\Tsim \to \infty }\J_\Tsim(\g) = \J(\g),%
\end{align*}
for all $\g \in \cS \subset \cK$.
Therefore, since the series of real numbers $\{\J_\Tsim(\g)\}_{\Tsim \in \N}$ converges to $\J(\g)$ and $\J(\g)$ is finite since $\g \in \cK$ by assumption, we can exploit the Cauchy convergence criterion to demonstrate that, for any $\acc > 0$ and $\cS \in \cK$, there exists a finite $\bar{\Tsim}\in \N$, possibly function of $\acc$ and $\cS$, such that for any $\Tsim > \bar{\Tsim}$, the bound~\eqref{eq:bound_sim} is achieved and the proof concludes.

\subsection{Proof of Lemma~\ref{lemma:estimation}}
\label{sec:estimation proof}

We note that~\cite[Lemma~1]{mimmo2021extremum} provides the same results claimed in Lemma~\ref{lemma:estimation}.
The only difference is that, in the mentioned reference, the objective function is assumed to be globally $\cC^3$.
However, since we assumed compactness of the set $\cS \subset \cK$ and since $\J$ and its gradient $\G$ are continuously differentiable~\cite{bu2019lqr} over the set of stabilizing gains $\cK$, 
we can repeat all the steps in~\cite[Lemma~1]{mimmo2021extremum} to get the proof over $\cS$.

\subsection{Proof of Lemma~\ref{lemma:averaged}}
\label{sec:proof_lemma_averaged}

Let us start by using the cost $\J$ to introduce the function $V_{\g}: \cK \to \R$ defined as 
\begin{align}
	V_{\g}(\ga) := \J(\ga) - \J(\gstar).
	\label{eq:V_g}
\end{align}
Being $\gstar$ the unique minimizer of $\J$~\cite{bu2019lqr}, we note that $V_{\g}$ is positive definite.
Now, given any $c > 0$, let us introduce 
\begin{subequations}\label{eq:levelset}
	\begin{align}
		\tlevset_c\! &:=\! \big\{(\fa,\ga) \in \R \times \cK \mid
		\notag\\
		&\hspace{.91cm} 
		\tfrac{1}{2}\norm{\fa \!-\! \Ja(\ga)}^2 \!+\! \J(\ga) \!-\! \J(\gstar) \!\leq\! c\big\}
		\\
		\levsetg_c \!&:=\! \big\{\ga \in \cK \mid \J(\ga) - \J(\gstar)\leq c \big\},
		\label{eq:levelsetg}
	\end{align}
\end{subequations}
namely, $\levsetg_c$ is the level set of $V_{\g}$ (cf.~\eqref{eq:V_g}), while $\tlevset$ would be the level set of $V$ (cf.~\eqref{eq:V_definition}) in the case in which $\mav = 1$.
Then, let $c_0 >0$ be the smallest number such that $(\fa\ud{0} - \Ja(\ga\ud{0}),\ga\ud{0}) \in \tlevset_{c_0}$ and use $\levsetg_{c_0}$ to define
\begin{align}
		\lipp_0 &:= 
		\max_{\ga \in \levsetg_{c_0}}\norm{\G(\ga)}.
		\label{eq:lipp_0}
\end{align}
We remark that~\cite[Corollary~3.7.1]{bu2019lqr} guarantees that, given any $c > 0$, the level set of the cost function $\J$, namely $\{\ga \in \R^{m \times n} \mid \J(\ga) \leq c\} \subset \R^{m \times n}$, is compact and, thus, so is $\levsetg_{c_0}$.
Hence, by continuity and differentiability of $\J$ and $\G$~\cite{bu2019lqr}, $\lipp_0$ is finite.
Now, by considering the compact set $\levsetg_{c_0}$ and $\delta \in (0,1]$, we recall that~\eqref{eq:ell_bound} (cf. Lemma~\ref{lemma:estimation}) ensures the existence of $\lippl > 0$ such that $\norm{\err(\gat)} \leq \lippl$ and that, for any $\acc > 0$, the result~\eqref{eq:bound_sim} (cf. (cf. Lemma~\ref{lemma:finite_time_sim})) ensures the existence of $\bar{\Tsim} > 0$ such that, for all $\Tsim > \bar{\Tsim}$, it holds $|\J(\g) - \Jt(\g)| \leq \acc$.
By exploiting these results, the definition of $\app$ (cf.~\eqref{eq:app}), and the triangle inequality, we write
\begin{align}
	\norm{\app(\g)} \leq \delta^2\lippl + \acc \dfrac{2}{\delta \Tper}\sum_{\initer=\iter+1}^{\iter + \Tper}\|\dr\ud\initer\|,\label{eq:bound_app_int}
\end{align}
for all $\g \in \levsetg_{c_0}$ and $\delta \in (0,1]$.
Now, to simplify the computations, we impose $\acc = \delta^3$.
We remark that, for all $\delta > 0$, this choice of $\acc$ is justified by Lemma~\ref{lemma:finite_time_sim} with a sufficiently large $\Tsim$.
In any case, this choice allows us to rewrite~\eqref{eq:bound_app_int} as 
\begin{align}
	\norm{\app(\g)} \leq \delta^2\lippa,\label{eq:bound_app}
\end{align}
for all $\g \in \levsetg_{c_0}$, where $\lippa :=\lippl  +2\sum_{\initer=\iter+1}^{\iter + \Tper}\|\dr\ud\initer\|/\Tper$.
Hence, by using~\eqref{eq:lipp_0},~\eqref{eq:bound_app}, and the triangle inequality, we get
\begin{align}
	\norm{\G(\ga) - \app(\ga)} \leq \lipp_0 + \delta^2\lippa,%
\end{align}
for all $\ga \in \levsetg_{c_0}$ and $\delta \in (0,1]$.
Thus, since $\cK$ is open~\cite[Lemma~IV.3]{bu2020topological}, for any $\tilde{c}_0 > c_0$, there exists $\bar{\ssz}_0 > 0$ such that 
\begin{align}
	\ga - \ssz\G(\ga) + \ssz\app(\ga) \in \levsetg_{\tilde{c}_0} \subset \cK,\label{eq:levset_c_1}
\end{align} 
for all $\ssz \in (0,\bar{\ssz}_0)$, $\delta \in (0,1]$, and $\ga \in \levsetg_{c_0}$.
We now invoke~\cite[Lemma~3.12]{bu2019lqr} to guarantee that the cost $\J$ is gradient dominated, i.e., there exists $\gd > 0$ such that
\begin{align}
	\J(\ga)  -  \J(\gstar) \leq \gd \norm{\G(\ga)}^2,\label{eq:gradient_dominance}
\end{align}
for all $\ga \in \cK$.
Now, we define 
	\begin{align}\label{eq:lipp}
		\lipp_1 &:= 
		\max_{\ga \in \levsetg_{\tilde{c}_0}}\norm{\G(\ga)}
		, \hspace{.15cm}
		\lipp_2  :=
		\max_{\ga \in \levsetg_{\tilde{c}_0}}\norm{\nabla G(\ga)}\!.
	\end{align}
Since also $\levsetg_{\tilde{c}_0}$ is compact~\cite[Corollary~3.7.1]{bu2019lqr} and recalling the continuity and differentiability of $\J$ and $\G$~\cite{bu2019lqr}, $\lipp_1$ and $\lipp_2$ are finite.
Next, we will use them to show that $\R\times \levsetg_{c_0}$ is forward-invariant for~\eqref{eq:averaged_system}.
To this end, assume that $\ga \in \levsetg_{c_0}$ and let us use an induction argument.
The increment $\Delta V_{\ga}(\ga)$ of $V_{\ga}(\gat)$ along trajectories of~\eqref{eq:averaged_system_gt} is given by
\begin{align}
	\Delta V_{\ga}(\gat)
	& := 
	\J(\ga - \ssz\G(\ga)+\ssz\app(\ga)) - \J(\ga)
	\notag\\
	&\stackrel{(a)}{\leq} 
	-\ssz\norm{\G(\ga)}^2 + \ssz\norm{\G(\ga)}\norm{\app(\ga)} 
	\notag\\
	&\hspace{.4cm}
	+\ssz^2\frac{\lipp_2}{2}\norm{\G(\ga) - \app(\ga)}^2
	\notag\\
	&\stackrel{(b)}{\leq} 
	-\ssz\left(1 - \ssz\lipp_2\right)\norm{\G(\ga)}^2
	\notag\\
	&\hspace{.4cm}
	+\ssz\norm{\G(\ga)}\norm{\app(\ga)} 
	\notag\\
	&\hspace{.4cm}
	+ \ssz^2\lipp_2\norm{\app(\ga)}^2
	, 
	\label{eq:deltaVav_intermediate}
\end{align}
where $(a)$ uses the Taylor expansion of $\J(\cdot)$ about $\ga$ evaluated at $\g - \ssz \G(\ga) + \ssz\app(\ga)$,~\eqref{eq:levset_c_1},~\eqref{eq:lipp}, and the Cauchy-Schwarz inequality, while $(b)$ rearranges the terms and uses $\frac{1}{2}\norm{\G(\ga) - \app(\ga)}^2 \leq \norm{\G(\ga)}^2 + \norm{\app(\ga)}^2$.
Let us arbitrarily fix $\eta \in (0,1)$ and define $\bar{\ssz}_2 := \min\{\bar{\ssz}_1,\frac{1 -\eta}{\lipp_2}\}$.
Then, for all $\ssz \in (0,\bar{\ssz}_2)$, we can bound~\eqref{eq:deltaVav_intermediate} as
\begin{align}
	\Delta V_{\ga}(\ga)  &\leq - \ssz\eta\norm{\G(\ga)}^2 +\ssz\norm{\G(\ga)}\norm{\app(\ga)} 
	\notag\\
	&\hspace{.4cm}
	+ \ssz^2\lipp_2\norm{\app(\ga)}^2
	\notag\\
	&\stackrel{(a)}{\leq}
	- \ssz\eta\norm{\G(\ga)}^2 
	+
	\ssz\lipp_0\delta^2\lippa^2
	\notag\\
	&\hspace{.4cm}
	+ \ssz^2\delta^4\lipp_2\lippa^2
	,\label{eq:deltaVav_4_intermediate}
\end{align}
where in $(a)$ we use~\eqref{eq:lipp_0} and~\eqref{eq:bound_app} to bound $\G(\ga)$ and $\norm{\app(\ga)}$ over the compact set $\levsetg_{c_0}$.
Now, in order to handle also the dynamics~\eqref{eq:averaged_system_ft}, let us introduce $\tfa \in \R$ defined as 
\begin{align}
	\tfa := \fa - \Ja(\ga),
\end{align}
which allows us to rewrite~\eqref{eq:averaged_system} as 
\begin{subequations}\label{eq:averaged_system_error}
	\begin{align}
		\tfatp &= (1-\ssz)\tfat + \pert(\gat)\label{eq:tf_dyn}
		\\
		\gatp &= \gat - \ssz\G(\gat) + \ssz\app(\gat),
	\end{align}
\end{subequations}
where $\pert: \R^{\dimu \times \dimx} \to \R$ is defined as 
\begin{align}\label{eq:pert}
	\pert(\g) := \Ja(\g)-\Ja(\g - \ssz\G(\g)+ \ssz\app(\g)).
\end{align}
Now, let us introduce $V_\f: \R \to \R$ defined as 
\begin{align}
	V_\f(\tfa) = \norm{\tfa}^2/2.
\end{align}
Hence, the increment $\Delta V_\f(\tfa) := V_\f((1-\ssz)\tfa + \pert(\ga)) - V_\f(\tfa)$ of $V_\f$ along the trajectories of~\eqref{eq:tf_dyn} reads as 
\begin{align}
	\Delta V_\f(\tfa) &= -\ssz(1 - \ssz/2)\norm{\tfa}^2 
	+ (1-\ssz)\tfat\pert(\ga) 
	\notag\\
	&\hspace{.4cm}+ 1/2\norm{\pert(\ga)}^2.
	\label{eq:Delta_V_tf}
\end{align}
Being the set $\cK$ open~\cite[Lemma~IV.3]{bu2020topological} and since $\drt$ is uniformly bounded for all $\iter$, there exists $\bar{\delta}_1 > 0$ such that $\ga + \delta\drt \in \cK$ for all $\ga \in \levsetg_{c_0}$, $\delta \in (0,\bar{\delta}_1)$, and $\iter \in \N$.
Hence, by exploiting the same arguments used to derive~\eqref{eq:lipp}, there exists $\lippp > 0$ such that 
\begin{align}
	\min_{\ga \in \tlevset_{c_0}}\norm{\nabla\Ja(\ga)} \leq \lippp,\label{eq:lipp_4}
\end{align}
for all $\delta \in (0,\bar{\delta}_1)$.
Thus, by using the definition of $\pert$ (cf.~\eqref{eq:pert}) and the triangle inequality, the bound~\eqref{eq:lipp_4} leads to
\begin{align}
	\norm{\pert(\ga)} \leq \ssz\lippp\norm{\G(\ga)} + \ssz\lippp\norm{\app(\ga)},\label{eq:bound_pert}
\end{align}
for all $\ga \in \Omega_{c_0}$ and $\delta \in (0,\bar{\delta}_1)$.
Hence, by using~\eqref{eq:bound_pert}, the Cauchy-Schwarz inequality, and the Young's inequality with parameter $2$, we can bound~\eqref{eq:Delta_V_tf} as 
\begin{align}
	\Delta V_\f(\tfa) &\leq
	-\ssz\big(1 - \tfrac{\ssz}{2}\big)\norm{\tfa}^2 
	+ 
	(1 -  \ssz)\ssz\lippp\!\!\norm{\tfa}\!\norm{\G(\ga)}
	\notag\\
	&\hspace{.4cm}
	+ \ssz^2\lippp^2\norm{\G(\ga)}^2
	\notag\\
	&\hspace{.4cm}
	+
	(1 -  \ssz)\ssz\lippp\norm{\tfat}\norm{\app(\ga)}
	\notag\\
	&\hspace{.4cm}
	+ 
	\ssz^2\lippp^2\norm{\app(\ga)}^2
	\notag\\
	&\stackrel{(a)}{\leq}\!\!
	-\ssz\big(1 - \tfrac{\ssz}{2}\big)\norm{\tfa}^2 
	+ 
	(1 -  \ssz)\ssz\lippp\!\!\norm{\tfa}\!\norm{\G(\ga)}
	\notag\\
	&\hspace{.4cm}
	+ \ssz^2\lippp^2\norm{\G(\ga)}^2
	\notag\\
	&\hspace{.4cm}
	+
	\delta^2(1 -  \ssz)\ssz\lippp\norm{\tfa}\lippl
	+ 
	\delta^4\ssz^2\lippp^2\lippl^2
	\notag\\
	&\stackrel{(b)}{\leq}
	-\ssz(1-\kappa/2 - \ssz/2)\norm{\tfa}^2 
	\notag\\
	&\hspace{.4cm}
	+ 
	(1 -  \ssz)\ssz\lippp\norm{\tfa}\norm{\G(\ga)}
	\notag\\
	&\hspace{.4cm}
	+ \ssz^2\lippp^2\norm{\G(\ga)}^2
	\notag\\
	&\hspace{.4cm}
	+ 
	\delta^4\left(\ssz^2 + \ssz/(2\kappa)(1 -  \ssz)^2\right)\lippp^2\lippl^2
	,\label{eq:DeltaV_tf_final}
\end{align}
where in $(a)$ we use~\eqref{eq:bound_app} to bound $\norm{\app(\ga)}$, while in $(b)$ we use the Young's inequality with an arbitrarily fixed parameter $\kappa \in (0,2)$ and write $\delta^2(1 -  \ssz)\ssz\lippp\norm{\tfa}\lippl \leq \ssz(\frac{\kappa}{2}\norm{\tfa}^2 + \frac{1}{2\kappa}\delta^4(1 -  \ssz)^2\lippp^2\lippl^2)$.
Now, let us compactly write all the terms due to the approximation error $\app$ by introducing $\tV: \R \times \R \times \R \to \R$ defined as
\begin{align}
	\hspace{-.25cm}
	\tV(\ssz,\delta,\mav) \! := \! \delta^2\lipp_0\lippa^2 \! + \! \delta^4(\ssz\lipp_2\lippa^2 \! + \! \tfrac{1}{\mav}(\ssz \!+ \! \tfrac{1}{2\kappa}\lippp^2\lippa^2)).\!\!
	\label{eq:tilde_V}
\end{align}
Then, let us $\V$ (cf.~\eqref{eq:V_definition}) and evaluate its increment $\Delta \V(\tfa,\ga) := \V((1-\ssz)\tfa + \pert(\ga),\ga - \ssz\G(\ga) + \ssz\app(\ga)) - \V(\tfa,\ga)$ along the trajectories of~\eqref{eq:averaged_system_error}.
By using~\eqref{eq:deltaVav_4_intermediate},~\eqref{eq:DeltaV_tf_final}, and the definition of $\tV$ (cf.~\eqref{eq:tilde_V}), we get 
\begin{align}
	&\Delta \V(\tfa,\ga)
	\notag\\
	&\leq -\ssz\begin{bmatrix}
		\tfa
		\\
		\G(\ga)
	\end{bmatrix}\T \cQ(\ssz,\mav)
	\begin{bmatrix}
		\tfa
		\\
		\G(\ga)
	\end{bmatrix} + \ssz\tV(\ssz,\delta,\mav),\label{eq:deltaV_Q}
\end{align}
where we introduced the matrix $\cQ(\ssz,\mav) \in \R^{2 \times 2}$ defined as 
\begin{align*}
	\cQ(\ssz,\mav) := \begin{bmatrix}
		\frac{1}{\mav}\big(1 - \frac{\kappa}{2} - \frac{\ssz}{2}\big)& -\frac{(1 - \ssz)\lippp}{2\mav}
		\\
		-\frac{(1 - \ssz)\lippp}{2\mav}&
		\eta - \ssz\frac{\lippp^2}{\mav} 
	\end{bmatrix}.
\end{align*}
Let us impose the positive definiteness of the top-left entry of $\cQ(\ssz,\mav)$.
To this end, let us arbitrarily fix $\nu \in (0,1-\kappa/2)$ and define $\bar{\ssz}_1 := \min\{\bar{\ssz}_0,\bar{\ssz}_2,2(1-\kappa/2-\nu)\}$.
Then, by Sylvester Criterion, for all $\ssz \in (0,\bar{\ssz}_1)$, it holds
\begin{align*}
	\cQ(\ssz,\mav) \succeq 
	\begin{bmatrix}
		\frac{\nu}{\mav}& -\frac{(1 - \ssz)\lippp}{2\mav}
		\\
		-\frac{(1 - \ssz)\lippp}{2\mav}& \eta - \ssz\frac{\lippp^2}{\mav}
	\end{bmatrix}.
\end{align*}
Now, let us impose the positive definiteness of $\cQ(\ssz,\mav)$.
To this end, we arbitrarily fix $\tilde{\eta} \in (0,\eta)$, $\tilde{\nu} \in (0,\nu)$, and define
\begin{align}
	\bar{\mav} := \max\left\{\frac{(1 + \bar{\ssz}_1^2)\lippp^2 + 4\bar{\ssz}_1\lippp^2(\nu - \tilde{\nu})}{4(\eta - \tilde{\eta})(\nu - \tilde{\nu})},1\right\}.
\end{align}
Then, we arbitrarily fix $\mav \ge \bar{\mav}$ and the Sylvester Criterion yelds
\begin{align*}
	\cQ(\ssz,\mav) \succeq \begin{bmatrix}
		\frac{\tilde{\nu}}{2\mav}& 0
		\\
		0& \tilde{\eta}
	\end{bmatrix},
\end{align*}
which allows us to further bound the right-hand side of~\eqref{eq:deltaV_Q} as 
\begin{align}
	\Delta \V(\tfa,\ga) 
	&\leq -\ssz (\tilde{\nu}/\mav\norm{\tfa}^2 + \tilde{\eta}\norm{\G(\ga)}^2)
	\notag\\
	&\hspace{.4cm}
	+ \tV(\bar{\ssz}_1,\delta,\mav).
	\label{eq:last_inequality_lyap}
\end{align}
The gradient dominance property of $\J$ (cf.~\eqref{eq:gradient_dominance}) leads to
\begin{align}
	& - \tilde{\nu}/\mav\norm{\tfa}^2 -\tilde{\eta}\norm{\G(\ga)}^2
	\notag\\
	&\leq
	- \tilde{\nu}/\mav\norm{\tfa}^2 - \tilde{\eta}/\gd\left(\J(\ga) - \J(\gstar)\right)
	\notag\\
	&
	\stackrel{(a)}{\leq} 
	-a V(\tfa,\ga),\label{eq:recovering_V}
\end{align}
where in $(a)$ we use the definition of $\V$ (cf.~\eqref{eq:V_definition}) and $a := \min\{2\tilde{\nu},\tilde{\eta}/\gd\}$.
Then, by using~\eqref{eq:recovering_V}, we further bound~\eqref{eq:last_inequality_lyap} as
\begin{align}
	\Delta \V(\tfa,\ga) \leq - \ssz a\V(\tfa,\ga) + \ssz\tV(\ssz,\delta,\mav).
	\label{eq:delta_V_leq_minus_V}
\end{align}
Now, without loss of generality, we assume $\rav \leq c_0$.
Indeed, one may always recover such a condition by using $\max\{\rav,c_0\}$ in place of $c_0$.
Then, let us define $\bar{\delta}_3 > 0$ as
\begin{align*}
	\bar{\delta}_3 :=  \min\left\{\sqrt{\frac{a\rav}{\delta^2\lipp_0\lippa^2 \! + \! \delta^4(\bar{\ssz}_1\lipp_2\lippa^2 \! + \! \tfrac{1}{\mav}(\bar{\ssz}_1 \!+ \! \tfrac{1}{2\kappa}\lippp^2\lippa^2))}},1\right\}.
\end{align*}
Then, for all $\delta \in (0,\bar{\delta}_1)$ with $\bar{\delta}_1 := \min\{\bar{\delta}_2,\bar{\delta}_3\}$, the definition of $\tV$ (cf.~\eqref{eq:tilde_V}) allows us to bound~\eqref{eq:delta_V_leq_minus_V} as
\begin{align}
	\Delta \V(\tfa,\ga)  &\leq - \ssz a\V(\tfa,\ga) + \ssz a\rav.
	\label{eq:deltaVav_4}
\end{align}
Although~\eqref{eq:deltaVav_4} seems to conclude the proof, we recall that it has been obtained by assuming $\ga \in \tlevset_{c_0}$.
In other words, since $(\tfa\ud0,\ga\ud0) \in \tlevset_{c_0}$ by definition of $c_0$, to conclude the proof we only need to prove that the set $\R \times \levsetg_{c_0}$ is forward-invariant for system~\eqref{eq:averaged_system_error}.
To this end, consider $(\tfat,\gat) \in \tlevset_{c_0}$ and, in light of the definition of $\V$ (cf.~\eqref{eq:V_definition}), we note that 
\begin{align}
	\J(\gatp) - \J(\gstar) 
	&\leq \V(\tfatp,\gatp) 
	\notag\\
	&\stackrel{(a)}{\leq} 
	\V(\tfat,\gat)
	\label{eq:V_invariance}\\
	&\stackrel{(b)}{\leq} 
	1/2\norm{\tfat}^2+\J(\gat) - \J(\gstar)	
	\stackrel{(c)}{\leq}
	c_0
	,\notag
\end{align}
where in $(a)$ we use the fact that the right-hand side of~\eqref{eq:deltaVav_4} is non-positive for all $(\tfat,\gat) \in \R \times \cK$ such that $\V(\tfat,\gat) \ge \rav$, in $(b)$ we use the fact that $\mav\ge 1$, while $(c)$ follows by the definition of $\tlevset_{c_0}$ (cf.~\eqref{eq:levelset}) and that $(\tfat,\gat) \in \tlevset_{c_0}$ by hypothesis.
The inequality~\eqref{eq:V_invariance} proves the desired invariance property of $\R \times \levsetg_{c_0}$ and the proof concludes.

\end{document}